\documentclass[a4paper,11pt]{amsart}

\usepackage{graphicx} 

\usepackage[left=2.7cm,right=2.7cm,top=3.5cm,bottom=3cm]{geometry}
\usepackage{amsthm,amssymb,amsmath,amsfonts,mathrsfs,amscd,verbatim}
\usepackage[utf8]{inputenc}
\usepackage[T1]{fontenc}
\usepackage[all,cmtip]{xy}
\usepackage{latexsym}
\usepackage{longtable}
\usepackage{mathtools}
\usepackage[pagebackref]{hyperref}
\usepackage{graphicx}
\usepackage{cite}
\usepackage{hyperref}
\usepackage{tikz-cd}
\usepackage{bm}

\numberwithin{equation}{section}
\mathtoolsset{showonlyrefs}

\newfont{\cyr}{wncyr10 scaled 1100}
\newfont{\cyrr}{wncyr9 scaled 1000}

\theoremstyle{plain}
\newtheorem{theorem}{Theorem}[section]
\newtheorem{proposition}[theorem]{Proposition}
\newtheorem{lemma}[theorem]{Lemma}
\newtheorem{corollary}[theorem]{Corollary}

\theoremstyle{definition}
\newtheorem{definition}[theorem]{Definition}

\theoremstyle{remark}
\newtheorem{remark}[theorem]{Remark}

\DeclareMathOperator{\GL}{GL}

\newcommand{\new}{\mathrm{new}}

\newcommand{\N}{\mathrm{N}}

\makeatletter
\newsavebox{\@brx}
\newcommand{\llangle}[1][]{\savebox{\@brx}{\(\m@th{#1\langle}\)}%
  \mathopen{\copy\@brx\kern-0.5\wd\@brx\usebox{\@brx}}}
\newcommand{\rrangle}[1][]{\savebox{\@brx}{\(\m@th{#1\rangle}\)}%
  \mathclose{\copy\@brx\kern-0.5\wd\@brx\usebox{\@brx}}}
\makeatother

\theoremstyle{plain}

\newtheorem{defn/thm}[thm]{Definition/Theorem}
\newtheorem{defn/prop}[thm]{Definition/Proposition}

\title{Sun's conjecture on a certain ternary quadratic form}
\author{Shilun Wang}
\address{The department of mathematics, Tongji University, Shanghai 200092, China}
\email{slwang@tongji.edu.cn}
\date{}

\begin{document}

\maketitle

\tableofcontents

\section{Introduction}

Let
\[
P(x,y,z)= x(5x+1) +\dfrac{y(5y+1)}{2}+\dfrac{z(5z+1)}{2}.
\]
Sun conjectured that every integer $n>1$ can be represented by $P(x,y,z)$ over $\mathbb{Z}$; see \cite[Conjecture 5.3(i)]{Sun2024}. The purpose of this paper is to establish an effective asymptotic version of this conjecture.

\begin{theorem}\label{thm:effective-main}
There exists an effectively computable positive integer $N_0$ such that every integer $n\geq N_0$ can be written as
\[
n=x(5x+1)+\dfrac{y(5y+1)}{2}+\dfrac{z(5z+1)}{2}
\]
for some $x,y,z\in\mathbb{Z}$.
\end{theorem}
\begin{remark}
We prove that every nonnegative integer $n$ for which $20n+2$ is not squarefree is represented by $P(x,y,z)$. See Theorem \ref{thm:nonsquarefree}. Thus the size restriction in our argument is needed only in the squarefree case.
\end{remark}

Representation problems of this kind are closely related to ternary quadratic forms, modular forms of half-integral weight, and central values of twisted $L$-functions. Ono and Soundararajan \cite{OnoSoundararajan1997} studied
Ramanujan's ternary quadratic form $x^2+y^2+10z^2$ and proved that it represents every eligible non-squarefree integer. Qin \cite{Qin2017} investigated related representation problems, including those for $x^2+7y^2+49z^2$.
Wei Lu and Qin \cite{LuQin2017} extended the Hecke-theoretic method to ternary quadratic forms whose genera consist of two classes in the same spinor genus under suitable coprimality hypotheses.

A further ingredient comes from identities relating representation numbers of ternary and quaternary quadratic forms.
Qin \cite{Qin2022} used such identities to give a new proof of an alternative form of Tunnell's theorem. More recently, Jun Lu and Qin \cite{LuQin2025} applied
this method to the family of elliptic curves
\[
E_n\colon \quad y^2=x^3+2n^3,
\]
obtaining formulas for $L(E_n,1)$ and criteria for its nonvanishing.
These identities provide a way to identify Shimura lifts by comparing Hecke eigenvalues.

Our argument begins with an elementary equivalence: an integer $n\geq0$ is represented by $P$ if and only if
$20n+2$ is represented by
\[
f(x,y,z)=x^2+y^2+25z^2.
\]
See Proposition \ref{prop:Sun-equivalence}. 

We determine the genus of $f$, which consists of two integral equivalence classes represented by $f$ and
\[
g(x,y,z)=2x^2+3y^2+5z^2+2yz+2xz.
\]
The difference of their theta series gives the cusp form
\[
\psi=\dfrac{\theta_f-\theta_g}{4}
=\sum_{m\geq1}b(m)q^m
\in S_{3/2}(100,\chi_{\mathrm{triv}}).
\]
Using the Hurwitz quaternion order, we establish explicit identities between representation numbers of $f$ and $g$ at prime squares and those of two quaternary quadratic forms $\alpha$ and $\beta$ at primes.
These identities allow us to identify $\text{Sh}_1(\psi)$ with
\[
\phi=\dfrac{\theta_\alpha-\theta_\beta}{4},
\]
a normalized newform of weight two and level $50$.

The genus-average formula and Gauss's three-square formula give
\[
r_f(t)=\dfrac{12}{5}H(4t)+\dfrac{16}{5}b(t)
\qquad
(t\equiv2\pmod{20}),
\]
where $H$ denotes the Hurwitz class number. For squarefree $t$ in this progression, $H(4t)=h(-4t)$. We treat the non-squarefree case by Hecke relations and Deligne's bound.

For the squarefree case, we use Banks's results on Landau--Siegel zeros and the effective Brauer--Siegel theorem \cite{BanksLS,BanksBS}, together with Waldspurger's formula and the subconvexity estimate of Blomer and Harcos
\cite{BlomerHarcos,BlomerHarcosAddendum}.
These results yield 
\[
h(-4t)\gg\frac{\sqrt t}{\log(4t)},
\qquad
|b(t)|\ll t^{15/32},
\]
with effectively computable implied constants. Consequently, the cuspidal contribution is smaller than the genus average beyond an effectively computable threshold. Combining the squarefree and non-squarefree cases proves
Theorem \ref{thm:effective-main}.

The paper is organized as follows. Section \ref{sec:genus} determines the genus of $f$ and establishes the required genus relation between $\alpha$ and $\beta$.
Section \ref{sec:theta} introduces the associated theta series and their Hecke actions.
Section \ref{sec:representation-identities} proves the representation identities using the Hurwitz quaternion order.
Section \ref{sec:shimura} identifies the Shimura lift of $\psi$.
Section \ref{sec:waldspurger} applies Waldspurger's formula and treats the relevant local square classes.
Section \ref{sec:application} proves the representation results and applies them to Sun's conjecture.

\medskip
\noindent\textit{Use of AI.}
ChatGPT was used for computations, editing, and
reorganization.

\section{Quadratic forms and their genera}\label{sec:genus}
For a positive definite quadratic form $Q \in \mathbb{Z}[x_1, \cdots, x_r]$, we can write 
\begin{equation}
Q(\bm{x})=\dfrac{1}{2} \bm{x}^\mathsf{T} A_Q \bm{x}, \quad A_Q=\begin{pmatrix}
\dfrac{\partial^2 Q}{\partial x_i\partial x_j}
\end{pmatrix}. 
\end{equation}
Thus $A_Q$ has integral entries and even diagonal. 
 
We denote $R_Q(n)=\{
v\in \mathbb{Z}^r \mid n=Q(v)
\}$ and $r_Q(n)=\#R_Q(n)$. 

Let \begin{equation}
O(Q)=\{
U \in \GL_r(\mathbb{Z}) \colon U^\mathsf{T} A_Q U=A_Q
\}. 
\end{equation}
We denote $o(Q)=\# O(Q)$. 

 Two integral quadratic forms $Q$ and $Q^\prime$ of the same rank are \emph{integrally
equivalent} if
\[
 A_{Q^\prime}=U^{\mathsf{T}}A_QU \quad \text{for some }U \in \GL_r(\mathbb{Z}).
\]
We write $Q\sim_{\mathbb{Z}}Q^\prime$ if they are integrally equivalent. 
For a prime $\ell$, they are \emph{locally equivalent at $\ell$} if 
\[
 A_{Q^\prime}=U^{\mathsf{T}}A_QU \quad \text{for some }U \in \GL_r(\mathbb{Z}_\ell).
\]
We write $Q\sim_{\mathbb{Z}_\ell}Q^\prime$ if they are locally equivalent at $\ell$. 

\begin{definition}[Genus of an integral quadratic form]\label{def:genus}
Two positive-definite integral quadratic forms $Q$, $Q^\prime$ are in the
\emph{same genus} if $Q\sim_{\mathbb{Z}_\ell}Q^\prime$ for every prime $\ell$. 
\end{definition}

We write
\[
 \text{gen}(Q)=\bigl\{[Q^\prime]_{\mathbb{Z}} \colon  Q^\prime\sim_{\mathbb{Z}_\ell}Q
                  \text{ for every prime }\ell\bigr\},
\]
where $[Q^\prime]_{\mathbb{Z}}$ denotes the integral equivalence classes of $Q^\prime$.

\subsection{The genus of the two ternary forms}

Let $Q(x,y,z)=ax^2+by^2+cz^2+ryz+sxz+txy$ be a positive-definite integral ternary quadratic form. We associate to $Q$ the quadratic lattice $L_Q=(\mathbb{Z}^3, Q(\mathbf{x}))$. 
We can also associate to $Q$ the symmetric bilinear form
\begin{equation}
B_Q(\mathbf{x},\mathbf{y}) =\dfrac{1}{2}\bigl(Q(\mathbf{x}+\mathbf{y}) -Q(\mathbf{x})-Q(\mathbf{y})\bigr).
\end{equation}
Given a $\mathbb{Z}$-basis $\{e_1,e_2,e_3\}$ of $L_Q$, the \emph{Gram matrix} of $L_Q$ with respect to this basis is
\[
G_{L_Q}=\bigl(B_Q(e_i,e_j)\bigr)_{1\leq i,j\leq 3}.
\]
With respect to the standard basis of $\mathbb{Z}^3$, the Gram matrix $G_{L_Q}$ of $L_Q$ is $\frac{1}{2}A_Q$. We then define $\det (L_Q) \colon \!\!=\det (G_{L_Q})=\dfrac{1}{8}\det( A_Q)$. 

Let
\begin{equation}
L_Q^\#=\{
x\in L_Q \otimes_\mathbb{Z}\mathbb{Q}\colon B_Q(x,L_Q)\subseteq\mathbb{Z}
\}
\end{equation}
be the dual lattice of $L_Q$.

Now we consider 
\begin{align*}
 f(x,y,z)&=x^2+y^2+25z^2,\\
 g(x,y,z)&=2x^2+3y^2+5z^2+2yz+2xz.
\end{align*}

\begin{proposition}[The two-class ternary genus]\label{prop:genusfg}
The genus of $f$ consists of exactly two integral equivalence classes,
represented by $f$ and $g$. Their automorphism numbers are $o(f)=16$, $o(g)=4$.
\end{proposition}

\begin{proof}
We denote by $L_f$ the associated quadratic lattice of $f$. Let $L_Q$ be the quadratic lattice associated to the positive-definite integral ternary quadratic form $Q$ in $\text{gen}(f)$. Hence for every prime $p$, there exists a matrix $U_p\in \GL_3(\mathbb{Z}_p)$ such that $A_Q=U_p^{T}A_fU_p$. Taking determinants, we obtain
\[
\det( A_Q)=(\det (U_p))^2\det( A_f).
\]
Since $\det (U_p)\in\mathbb{Z}_p^\times$, it follows that $v_p(\det(A_Q))=v_p(\det( A_f))$ for every prime $p$. As both determinants are positive integers, they must be equal.  Hence $\det(A_Q)= \det( A_f)=200$. 
We have
\[
\det(L_Q)=\det( L_f)=25.
\]
Consider the quadratic lattice 
\[
M=L_Q+5L_Q^\#.
\]

Since $Q\in \text{gen}(f)$, we have $(L_Q)_p\simeq(L_f)_p$. Hence we can choose an orthogonal basis $e_1$, $e_2$, $e_3$ of $(L_Q)_p$ with Gram matrix $\text{diag}(1,1,25)$.

For $p\neq 5$, since $\det (L_Q)=25$ is a $p$-adic unit, we have $(L_Q^\#)_p=(L_Q)_p$. See more details in \cite[\S 82G]{OMeara}. 
So $M_p=(L_Q)_p+5(L_Q^\#)_p=(L_Q)_p$. Thus $G_{M_p} \in \GL_3(\mathbb{Z}_p)$. 

For $p=5$,  we have 
\[
(L_Q^\#)_5=\mathbb{Z}_5e_1\oplus\mathbb{Z}_5e_2 \oplus\mathbb {Z}_5\dfrac{e_3}{25}.
\]
Thus
\[
M_5 =\mathbb{Z}_5e_1\oplus \mathbb{Z}_5e_2  \oplus \mathbb{Z}_5\dfrac{e_3}{5}
\]
Since
\[
B_Q\left(\frac{e_3}{5},\frac{e_3}{5}\right) =\frac{1}{25}B_Q(e_3,e_3)=1
\]
and the basis vectors are pairwise orthogonal, the Gram matrix of $M_5$ with respect to $e_1$, $e_2$, $e_3/5$ is the identity matrix. Thus $G_{M_5}\in \GL_3(\mathbb{Z}_5)$.

Since $G_{M_p}\in \GL_3(\mathbb{Z}_p)$ for every prime $p$, the Gram matrix $G_M$ has entries in $\mathbb{Z}$ and $\det(G_M) \in \mathbb{Z}^\times=\{\pm1\}$. As $M$ is positive definite, we conclude that $\det(M)=1$. 
Since $M$ is a positive-definite integral lattice of rank three with $\det(M)=1$, the classification of such low-dimensional lattices implies that $M\simeq I_3$, where $I_3=(\mathbb{Z}^3, x^2+y^2+z^2)$. See more details in \cite[Theorem 1 and Table 1]{ConwaySloane1982}.
Since $\det (L_Q)=25$ and $\det(M)=1$, the determinant formula gives $[M: L_Q]=5$. Hence $L_Q$ is isometric to an index-five sublattice $N$ of $I_3$.
As $\mathbb{Z}^3/N$ is a group of order $5$, we have $\mathbb{Z}^3/N\simeq\mathbb{F}_5$. Thus $N$ is the kernel of a nonzero homomorphism $\mathbb{Z}^3\to\mathbb{F}_5$. Hence $N$ is of the form
\[
L_{\bm{a}} =\left\{
(x,y,z)\in\mathbb{Z}^3\colon a_1x+a_2y+a_3z\equiv0\pmod 5
\right\},
\]
for some $[\bm{a}]=[a_1:a_2:a_3] \in\mathbb{P}^2(\mathbb{F}_5)$.

We claim that such $L_{\bm{a}}$ are only given by $a_1^2+a_2^2+a_3^2\in\{1,4\}\pmod5$.

Since $[I_3:L_{\bm{a}}]=5$, we have $(L_{\bm{a}})_p=(I_3)_p\simeq(L_f)_p$ for every prime $p\ne5$. Hence it suffices to consider $p=5$. We set $u\colon \!\!=a_1^2+a_2^2+a_3^2$ and suppose that $u \not\equiv 0\pmod5$.
We choose a lift $w\in\mathbb{Z}_5^3$ of $\mathbf a$. We put $\widetilde{u}=B(w,w)\in\mathbb{Z}_5^\times$, where $B$ is the standard bilinear form on $I_3$. Hence $\widetilde{u}\equiv  a_1^2+a_2^2+a_3^2 =u \not\equiv 0\pmod 5$, which means that $\widetilde u$ is a $5$-adic unit. 

For any $x\in\mathbb{Z}_5^3$, we may write
\[
x=\dfrac{B(x,w)}{\widetilde{u}}w +\left(
x-\dfrac{B(x,w)}{\widetilde{u}}w
\right),
\]
where the first coefficient belongs to $\mathbb{Z}_5$. We can show that the  second summand $x-\dfrac{B(x,w)}{\widetilde{u}}w$ lies in $\{x\in\mathbb{Z}_5^3 \colon B(x,w)=0\}$. Consequently, there is a decomposition 
\begin{equation}\label{decomposition}
\mathbb{Z}_5^3=\mathbb Z_5w\oplus w^\perp, 
\end{equation}
where $w^\perp\colon \!\!= \{x\in\mathbb{Z}_5^3 \colon B(x,w)=0\}$. 

By the defining congruence of
$L_{\bm{a}}$, we have $(L_{\bm{a}})_5 =\{
x\in\mathbb{Z}_5^3 \colon B(x,w)\in 5\mathbb{Z}_5 
\}$.
We write $x=\lambda w+y$ with $y\in w^\perp$ as in Equation \eqref{decomposition}.  We get $B(x,w)=\lambda\widetilde{u}$. 
Since $\widetilde{u}\in\mathbb{Z}_5^\times$, the condition $B(x,w)\in 5\mathbb Z_5$ is equivalent to $\lambda\in 5\mathbb{Z}_5$.
Hence
\[
(L_{\bm{a}})_5 =5\mathbb{Z}_5 w  \oplus w^\perp.
\]
Obviously $w^\perp$ is a $\mathbb{Z}_5$-lattice of rank two. So we choose a $\mathbb{Z}_5$-basis $(v_1,v_2)$. Then the Gram matrix with respect to $(w,v_1,v_2)$ is
\[
G_{I_3}=
\begin{pmatrix}
\widetilde{u}&0&0\\
0&B(v_1,v_1)&B(v_1,v_2)\\
0&B(v_2,v_1)&B(v_2,v_2)
\end{pmatrix}.
\]
Since $\det I_3=1$, we obtain
\[
\det(w^\perp)\equiv \widetilde{u}^{-1} \pmod{ (\mathbb Z_5^\times)^2}.
\] 
By the classification of quadratic lattices over $\mathbb{Z}_5$ via Jordan decompositions (see \cite[Chapter 15, Section 7]{ConwaySloane1999}), the ranks and determinant square classes of the Jordan components determine their
isometry class. Consequently,
\[
(L_{\bm{a}})_5\simeq (L_f)_5 
\quad\Longleftrightarrow\quad
\widetilde{u}\in(\mathbb{Z}_5^\times)^2.
\]
By Hensel's lemma, we have $\widetilde u\in(\mathbb{Z}_5^\times)^2 \Longleftrightarrow u\in(\mathbb{F}_5^\times)^2=\{1,4\}$.

Suppose that $u=0$ in $\mathbb{F}_5$. Then $\bm{a}$ is a nonzero isotropic vector in $\mathbb{F}_5^3$. The restriction of the standard bilinear form $B$ to $\bm{a}^\perp=\{
x\in\mathbb{F}_5^3\colon\bm{a}\cdot x=0
\}$ has radical $\mathbb F_5\mathbf a$ and hence has rank one. Since a basis of $(L_{\bm{a}})_5$ may be chosen whose first two vectors reduce to a basis of $\bm{a}^\perp$ and whose third vector lies in
$5\mathbb{Z}_5^3$, its Gram matrix has rank one modulo $5$. On the other hand, the Gram matrix of $(L_f)_5$ has rank two modulo $5$. Thus $(L_{\bm{a}})_5 \not\simeq(L_f)_5$. 

The integral orthogonal group $O(I_3)$ consists of signed permutation matrices, since
the only integral vectors of norm one in $I_3$ are $\pm (1,0,0),\pm (0,1,0), \pm (0,0,1)$. 

We now classify the projective points $[\bm{a}]\in\mathbb{P}^2(\mathbb{F}_5)$ satisfying $a_1^2+a_2^2+a_3^2\in\{1,4\}\pmod5$. Since the squares in $\mathbb{F}_5$ are $0,1,4$, such a point has either exactly one or three nonzero coordinates. 

In the first case, it is equivalent under coordinate permutations to $[0:0:1]$. In the second case, the coordinate squares are, up to permutation, $(1,1,4)$ or
$(1,4,4)$. Multiplying by a nonzero scalar and changing signs shows that every such point lies in the orbit of $[1:1:2]$.
Thus the admissible points form exactly two orbits under $O(I_3)$, represented by $[0:0:1]$ and $[1:1:2]$. Points in the same orbit define isometric lattices, since $O(I_3)$ acts by
isometries on $I_3$. 

For the first representative, let
\[
v_1=(1,0,0),\qquad
v_2=(0,1,0),\qquad
v_3=(0,0,5).
\]
These vectors belong to $L_{(0,0,1)}$.  Their coordinate matrix
\[
P=
\begin{pmatrix}
1&0&0\\
0&1&0\\
0&0&5
\end{pmatrix}
\]
has determinant $5$. Hence they form a $\mathbb{Z}$-basis of $L_{(0,0,1)}$.
Its Gram matrix is
\[
P^TP=
\begin{pmatrix}
1&0&0\\
0&1&0\\
0&0&25
\end{pmatrix}.
\]
Consequently, $L_{(0,0,1)}\simeq L_f$. 

For the second representative, let
\[
v_1=(-1,1,0),\qquad
v_2=(-1,-1,1),\qquad
v_3=(0,1,2).
\]
These vectors belong to $L_{(1,1,2)}$.  Their coordinate matrix
\[
P=
\begin{pmatrix}
-1&-1&0\\
1&-1&1\\
0&1&2
\end{pmatrix}
\]
has determinant $5$. Hence they form a $\mathbb{Z}$-basis of $L_{(1,1,2)}$.
Its Gram matrix is
\[
P^TP=
\begin{pmatrix}
2&0&1\\
0&3&1\\
1&1&5
\end{pmatrix}.
\]
Consequently, $L_{(1,1,2)}\simeq L_g$.

We know that $f$ represents $1$ precisely at
$(\pm1,0,0)$ and $(0,\pm1,0)$. Hence $r_f(1)=4$. On the other hand,
\[
g(x,y,z)=2\left(x+\dfrac {z}{2}\right)^2 +3\left(y+\dfrac{z}{3}\right)^2+\dfrac{25}{6}z^2.
\]
\begin{enumerate}
\item If $z\ne0$, then $g(x,y,z)\ge25/6>1$. 
\item If $z=0$, then $g(x,y,0)=2x^2+3y^2$ cannot equal $1$. 
\end{enumerate}
Hence $r_g(1)=0$.
Since representation numbers are invariant under integral equivalence, $f$ and $g$
are not integrally equivalent. Consequently, $\text{gen}(f) =\{[f],[g]\}$.

We next compute their automorphism numbers. 

For $f$, the vectors of norm one are precisely $(\pm1,0,0)$ and $(0,\pm1,0)$. Every integral automorphism acts on their span by a signed permutation matrix and sends $(0,0,1)$ to
$(0,0,\pm 1)$. Hence $o(f)=2!\cdot 2^2\cdot 2=16$. 

For $g$, the vectors of norm two are exactly $(\pm 1, 0, 0)$ and those of norm three are exactly
$(0, \pm 1, 0)$. 
Let $U\in O(g)$. Then
\[
Ue_1=\varepsilon e_1,\qquad
Ue_2=\delta e_2,
\qquad \varepsilon,\delta\in\{\pm1\}.
\]
Writing $Ue_3=(u,v,w)^T$, we have
\[
U=
\begin{pmatrix}
\varepsilon&0&u\\
0&\delta&v\\
0&0&w
\end{pmatrix},
\qquad w=\pm1,
\]
where the latter follows from
$\det U=\pm1$.

Since $U$ preserves the Gram matrix $G_{L_g}$, we obtain $2u+w=\varepsilon$, $3v+w=\delta$. As $\delta$, $w\in\{\pm 1\}$ and $v\in\mathbb{Z}$, the second equality forces $v=0$, $w=\delta$. Hence $u=\frac{\varepsilon-\delta}{2}$. 

Conversely, for every $\varepsilon$, $\delta\in\{\pm1\}$, the matrix
\[
U_{\varepsilon,\delta}=
\begin{pmatrix}
\varepsilon&0&(\varepsilon-\delta)/2\\
0&\delta&0\\
0&0&\delta
\end{pmatrix}
\]
satisfies
\[
U_{\varepsilon,\delta}^T
G_{L_g}U_{\varepsilon,\delta}=G_{L_g}.
\]
Therefore $o(g)=4$. 
\end{proof}

\subsection{The genus of the two quaternary forms}
Let
\begin{align*}
 \alpha(x,y,z,w)&=x^2+y^2+25z^2+13w^2+xw+yw+25zw,\\
 \beta(x,y,z,w)&=3x^2+3y^2+5z^2+5w^2+3xz-xw-yz-3yw.
\end{align*}

\begin{proposition}\label{prop:genusab}
The positive-definite quaternary quadratic forms $\alpha$ and $\beta$ belong to the same genus, but represent distinct equivalence classes.
\end{proposition}
\begin{proof}
Consider the rational matrix
\[
U_{\alpha\beta}=\dfrac{1}{5}
\begin{pmatrix}
-8&2&-1&2\\
2&7&-1&2\\
0&-1&1&3\\
1&1&2&-4
\end{pmatrix}.
\]
Direct calculation shows that $\det(U_{\alpha\beta})=1$ and $U_{\alpha\beta}^\mathsf{T}A_\alpha U_{\alpha\beta}=A_\beta$. Since $U_{\alpha\beta} \in\GL_4(\mathbb{Z}_\ell)$
whenever $\ell\neq5$, this proves the desired local equivalence at every prime $\ell\neq 5$. It remains to consider $\ell=5$.
We can rewrite $\alpha$ and $\beta$ as 
\begin{equation}
\alpha=\left(x+\dfrac{w}{2}\right)^2 +\left(y+\dfrac{w}{2}\right)^2+25\left(z+\dfrac{w}{2}\right)^2 +25\left(\dfrac{w}{2}\right)^2
\end{equation}
and
\begin{equation}
\beta=3\left(x+\dfrac z2-\dfrac w6\right)^2+3\left(y-\dfrac z6-\dfrac w2\right)^2 +\dfrac{25}{6}(z^2+w^2).
\end{equation}
Consider the matrix
\[
P_\alpha=
\begin{pmatrix}
1&0&0&\frac12\\
0&1&0&\frac12\\
0&0&1&\frac12\\
0&0&0&\frac12
\end{pmatrix}
\in \GL_4(\mathbb{ Z}_5).
\]
By the above identity, the change of variables
\[
(X_1,X_2,X_3,X_4)^{\mathrm T}
=P_\alpha(x,y,z,w)^{\mathrm T}
\]
transforms $\alpha$ into
\[
X_1^2+X_2^2+25X_3^2+25X_4^2.
\]

Consider the matrix
\[
U_\beta=
\begin{pmatrix}
1&0&-\dfrac{1}{2}&\dfrac{1}{6}\\
0&1&\dfrac{1}{6}&\dfrac{1}{2}\\
0&0&1&0\\
0&0&0&1
\end{pmatrix}
\in \GL_4(\mathbb{Z}_5).
\]
Then
\[
U_\beta^\mathsf{T} A_\beta U_\beta
=
\begin{pmatrix}
6&0&0&0\\
0&6&0&0\\
0&0&25/3&0\\
0&0&0&25/3
\end{pmatrix}.
\]

Choose $v,t\in\mathbb{Z}_5^\times$ such that
\[
v^2=-\dfrac{2}{3},\qquad t^2=\dfrac{1}{6}.
\]
Such elements exist by Hensel's lemma. Consider the matrix
\[
P_\beta=
\begin{pmatrix}
1&-v&0&0\\
v&1&0&0\\
0&0&t^{-1}&0\\
0&0&0&t^{-1}
\end{pmatrix}
\in\mathrm{GL}_4(\mathbb{Z}_5).
\]

Since $1+v^2=1/3$ and $t^{-2}=6$, we have
\[
P_\beta^\mathsf{T}(U_\beta^\mathsf{T} A_\beta U_\beta) P_\beta=
\begin{pmatrix}
2&0&0&0\\
0&2&0&0\\
0&0&50&0\\
0&0&0&50
\end{pmatrix}.
\]
Consequently, $\beta$ is equivalent to $X_1^2+X_2^2+25X_3^2+25X_4^2$ over $\mathbb{Z}_5$. Hence $\alpha$ and $\beta$ are in the same genus. 

We now show that $\alpha$ and $\beta$ are different classes in the genus, which means that they are not equivalent over $\mathbb{Z}$.
From the above expression for $\alpha$, we obtain $r_\alpha(1)=4$, with the four representations given by $(\pm1,0,0,0)$, $(0,\pm1,0,0)$.
On the other hand, suppose that
$\beta(x,y,z,w)=1$ for some
$x,y,z,w\in\mathbb Z$.
Since
\[
\beta=
3\left(x+\dfrac z2-\dfrac w6\right)^2
+3\left(y-\frac z6-\frac w2\right)^2
+\frac{25}{6}(z^2+w^2),
\]
we necessarily have $z=w=0$.
It follows that
\[
3(x^2+y^2)=1,
\]
which is impossible for integers $x,y$.
Thus
\[
r_\beta(1)=0.
\]

Since integrally equivalent quadratic forms
have identical representation numbers, the
inequality
\[
r_\alpha(1)=4\neq0=r_\beta(1)
\]
implies that $\alpha$ and $\beta$ represent two distinct integral equivalence classes in the same genus.

\end{proof}

\begin{remark}\label{rem:quaternary-genus}
Proposition \ref{prop:genusab} only establishes that $\alpha$ and $\beta$ belong to the same genus and represent different integral equivalence classes. In fact, the genus of $\alpha$ consists of exactly four integral equivalence classes. However, we do not determine the full genus here. 
\end{remark}

\section{The theta series and modular forms}\label{sec:theta}
For a positive definite integral form $Q$, its level $N_Q$ is the least positive integer for which $N_QA_Q^{-1}$ is integral with even diagonal. 
Its theta character is 
\begin{equation}
\chi_Q(\cdot)=
\begin{cases}
\left( \dfrac{2\det A_Q}{\cdot} \right),   \quad r \text{ is odd}. \\
\left( \dfrac{(-1)^{r/2}\det A_Q}{\cdot} \right) ,  \quad r \text{ is even}  .
\end{cases}
\end{equation}
 
 We define the following theta series
\begin{equation}
\theta_Q(\tau)=\sum\limits_{n\geq 0} r_Q(n)q^n, \quad q=\mathrm{e}^{2\pi i \tau}.
\end{equation}
We know that $\theta_Q(\tau)\in M_{r/2}(N_Q, \chi_Q)$. See \cite[Theorem 10.1]{WangPei2012}.

\begin{lemma}
\begin{enumerate}
\item The theta series $\theta_f$ and $\theta_g$ are in $M_{3/2}(100, \chi_\text{triv})$. 
\item The theta series $\theta_\alpha$ and $\theta_\beta$ are in $M_2(50, \chi_\text{triv})$. 
\end{enumerate}
Here 
\begin{equation}
\chi_\text{triv}(d)= 
\begin{cases}
1 , \quad \text{if } \gcd(d, 10)=1.\\
0, \quad \text{if } \gcd(d, 10)>1.
\end{cases}
\end{equation}
\end{lemma}

For $F=\sum \limits_{n=0}^\infty c(n)q^n \in M_{3/2}(100, \chi_\text{triv})$ and primes $p\nmid 10$, the Hecke operator $T_{p^2}$ acts by  
\begin{equation}
T_{p^2}(F)=\sum \limits _{n=0}^\infty (c(p^2n)+ \left (
\dfrac{-n}{p} \right)
c(n)+ p c(n/p^2) )q^n
\end{equation}
where $c(n/p^2)=0$ if $p^2\nmid n$.

\subsection{The ternary theta series}
 We define 
\begin{equation}
\Theta=\dfrac{\theta_f/o(f)+\theta_g/o(g)}{1/o(f)+1/o(g)}=\dfrac{\theta_f+4\theta_g}{5}, 
\end{equation}
$\vartheta=\theta_f-\theta_g$ and $\psi=\dfrac{1}{4}\vartheta=\sum\limits_{n\geq 0} b(n)q^n$.  From \cite{Schiemann1997}, we know that $\vartheta \neq 0$,  since $f$ and $g$ are different classes in the genus.  Hence $\psi\neq 0$. 

\begin{lemma}
We have $\Theta \in E_{3/2}(100, \chi_\text{triv})$ and $\vartheta\in S_{3/2}(100, \chi_\text{triv})$. Furthermore, $b(1)=1$. 
\end{lemma}
\begin{proof}
See \cite[Theorem 10.3 and p.365]{WangPei2012} and \cite[p.19]{LuQin2017}. 
We know that $r_f(1)=4$ and $r_g(1)=0$. So $b(1)=\dfrac{1}{4}(r_f(1)-r_g(1))=1$. 
\end{proof}

Let $M(f)$ denote the $\mathbb{C}$-linear space spanned by $\theta_f$ and $\theta_g$. 
We define 
\begin{enumerate}
\item $E(f)=M(f)\cap E_{3/2}(100, \chi_\text{triv})$. 
\item $S(f)=M(f)\cap S_{3/2}(100, \chi_\text{triv})$. 
\end{enumerate}
Hence $E(f)=\mathbb{C}\Theta$ and $S(f)=\mathbb{C}\vartheta$. 

\begin{proposition}\label{prop:hecke}
For every prime $p\nmid10$, the two-dimensional space $M(f)$ is stable under $T_{p^2}$ and
\[
 T_{p^2}\Theta=(p+1)\Theta,\qquad
 T_{p^2}\psi=\lambda_p\psi
\]
for some $\lambda_p\in\mathbb{C}$.
\end{proposition}
\begin{proof}
By Eichler's commutation relation, cf. \cite[Hilfssatz 1]{SchulzePillot1984}, the space $M(f)$ is stable under the Hecke operators $T_{p^2}$ for every prime $p\nmid 10$. Thus $E(f)$ and $S(f)$ are stable under $T_{p^2}$. 
Hence $\Theta$ and $\vartheta$ are eigenforms of $T_{p^2}$ for any prime $p\nmid 10$. For each prime $p\nmid 10$, we let $\lambda_p$ be the eigenvalue of $\psi$ with respect to the Hecke operator $T_{p^2}$. 
We refer to \cite[Satz 1]{SchulzePillot1984} and \cite[Theorem 10.2]{WangPei2012} for more details about the eigenvalue of $\Theta$. 
\end{proof}

\subsection{The quaternary theta series}

We define $\phi=\dfrac{1}{4}(\theta_\alpha-\theta_\beta)=\sum \limits_{n\geq 0}a(n)q^n$, which is in $S_2(50,\chi_\text{triv})$ from \cite[p. 365]{WangPei2012}. We also have $a(1)=1$. 
 
\begin{proposition}\label{prop:dim50}
We have $\dim_\mathbb{C} S_2(50, \chi_\text{triv})=2$ and $S_2(50, \chi_\text{triv})=S_2^{\new}(50, \chi_\text{triv})$.
\end{proposition}
\begin{proof}
Recall the modular-curve genus formula
\begin{equation}\label{eq:genusX0}
 g(X_0(N))=1+\frac{\mu(N)}{12}-\frac{e_2(N)}4
                  -\frac{e_3(N)}3-\frac{c(N)}2,
\end{equation}
where $\mu(N)=[\mathrm{SL}_2(\mathbb{Z}):\Gamma_0(N)]$, $e_2,e_3$ count elliptic points, and $c(N)$ counts cusps. See \cite[Chapter 3]{DiamondShurman}. 

Since
\[
 \mu(50)=90,\quad e_2(50)=2,\quad e_3(50)=0,  \quad c(50)=\sum_{d\mid50}\varphi(\gcd(d,50/d))=12,
\]
we have $\dim S_2(50, \chi_\text{triv})=g(X_0(50))=2$. Every proper divisor $M$ of $50$ is one of $1,2,5,10,25$, and the same formula yields $g(X_0(M))=0$ for each. Hence no oldform enters from a proper divisor and the Atkin--Lehner--Li decomposition proves the assertion \cite[Chapter 5]{DiamondShurman}.
\end{proof}

\subsection{The unary theta subspace}

For a level $N$ and a Dirichlet character $\chi$, let $U(N,\chi)$ denote the subspace of $S_{3/2}(N,\chi)$ spanned by unary theta series. Following \cite[Section 2]{Purkait2013}, this space has a basis consisting
of the forms \begin{equation} 
u_{\eta,t}(z)=\sum_{m\geq1}\eta(m)m q^{tm^2},
\end{equation}
where $t$ is a positive integer, $\eta$ is a primitive odd Dirichlet character of conductor $r_\eta$ and
\begin{equation}
4r_\eta^2t\mid N,
\qquad
\chi=
\left(\frac{-4t}{\cdot}\right)\eta.
\end{equation}

\begin{lemma}\label{lem:no-unary}
The unary theta subspace $U(100,\chi_\text{triv})=\{0\}$. 
\end{lemma}
\begin{proof}
Suppose that a unary theta series $u_{\eta,t}$ occurs in $U(100,\chi_\text{triv})$.  Then $4r_\eta^2t\mid100$, i.e., $r_\eta^2t\mid25$. Since $\eta$ is a primitive odd Dirichlet character, its conductor cannot be $1$. The divisibility condition therefore forces $r_\eta=5$ and $t=1$. The character condition becomes $\chi_\text{triv}=\left(\frac{-4}{\cdot}\right)\eta$. However, evaluating at $11$ gives $\left(\frac{-4}{11}\right)\eta(11)=(-1)\eta(1)=-1$, whereas the trivial character modulo $100$ takes the value $1$ at $11$. This contradiction shows that no such unary theta series exists. Therefore $U(100,\chi_\text{triv})=\{0\}$. 
\end{proof}

\section{Representation identities}\label{sec:representation-identities}
We establish explicit identities between the representation numbers of certain ternary quadratic forms at prime squares and those of quaternary quadratic forms at primes.
Our argument follows the method of Qin \cite{Qin2022} and Lu--Qin \cite{LuQin2025}. 

\subsection{Factorization in the Gaussian integers}
\label{subsec:gaussian-factorization}

Let \begin{equation}
\mathbb{Z}[i]=\{a+bi:a,b\in\mathbb{Z} \}, \qquad \N(a+bi)=a^2+b^2.
\end{equation}

\begin{lemma}\label{lem:gaussian-factorization}
Let $A,B$ be coprime positive integers and let $u\in\mathbb Z[i]$ satisfy $\N(u)=AB$. Then the system
\begin{equation}\label{eq:gaussian-system}
\xi\zeta=u,\qquad \N(\xi)=A,\qquad \N(\zeta)=B
\end{equation}
has exactly four solutions
$(\xi,\zeta)\in\mathbb Z[i]\times \mathbb{Z}[i]$.
\end{lemma}

\subsection{A quaternionic norm map}\label{subsec:quaternion-norm}
Let $\mathbb{H}$ be Hamilton's quaternion algebra over $\mathbb{Q}$, which we can write as $\mathbb{Q}+\mathbb{Q}i+\mathbb{Q}j+\mathbb{Q}k$,  where $i^2=j^2=-1$, $k=ij$ and $ij=-ji$. 
For $\rho=a+bi+cj+dk\in\mathbb{H}$, we write $\overline{\rho}=a-bi-cj-dk$ and then 
\begin{equation}
\N(\rho)=\rho \overline{\rho}=a^2+b^2+c^2+d^2.
\end{equation}
Let $\mathcal H$ be the Hurwitz order given by 
\begin{equation}
\mathcal{H}=\mathbb{Z}+\mathbb{Z}i+\mathbb{Z}j+\mathbb{Z}h, \qquad h=\frac{1+i+j+k}{2}.
\end{equation}
We define $T(\rho)=\rho i\overline{\rho}$. More precisely, we can write 
\begin{equation}
T(\rho)=Xi+Yj+Zk,
\end{equation}
where $X=(a^2+b^2-c^2-d^2)$, $Y=2(ad+bc)$ and $Z=2(bd-ac)$. We know that $(X,Y,Z)\in \mathbb{Z}^3$ for every $\rho\in\mathcal{H}$. Taking the norm, we have $\N(T(\rho))=\N(\rho)^2$, i.e., 
\begin{equation}\label{eq:hopf-identity}
X^2+Y^2+Z^2
=
(a^2+b^2+c^2+d^2)^2.
\end{equation}
 Consider the map $\mathcal{T}\colon \mathcal{H}\longrightarrow\mathbb{Z}^3$ defined by $\rho \longmapsto (X,Y,Z)$. For odd primes $p$, we set 
\[
 \mathcal{F}_p(X,Y,Z)=\{\rho\in\mathcal{H}\colon \N(\rho)=p, \mathcal{T}(\rho)=(X, Y, Z)\}.
\]

\begin{lemma}\label{lem:primitive-fibres}
 Let $p$ be an odd prime. Suppose $X^2+Y^2+Z^2=p^2$ and $\gcd(X,Y,Z)=1$. Then $\#\mathcal{F}_p(X,Y,Z)=4$.
\end{lemma}

\begin{proof}
Since $p$ is odd, we have $X^2+Y^2+Z^2\equiv1\pmod4$. Hence exactly one of $X,Y,Z$ is odd. Conjugation by $h$ cyclically permutes $(X,Y,Z)$. Thus it suffices to treat the case that $X$ is odd and $Y$, $Z$ are even.

If $\rho=a+bi+cj+dk \in \mathcal{H}$ is in $\mathcal{F}_p(X,Y,Z)$, then $a$, $b$, $c$, $d\in\mathbb{Z}$. Otherwise we write 
\[
a=\dfrac{a_0}{2},\quad
b=\dfrac{b_0}{2},\quad
c=\dfrac{c_0}{2},\quad
d=\dfrac{d_0}{2},
\]
where all four numerators are odd.
Since every odd square is congruent to $1$ modulo $8$, we have
\begin{equation}
a_0^2+b_0^2-c_0^2-d_0^2 \equiv 0 \pmod 8.
\end{equation}
Therefore $X$ is even, a contradiction. 

Set
\begin{equation}\label{eq:ABu}
A=\dfrac{p+X}{2},\qquad
B=\dfrac{p-X}{2},\qquad
u=\dfrac{Y+iZ}{2}.
\end{equation}
Since $\gcd(X, Y, Z)=1$, we have $X\neq\pm p$. Thus $\gcd(A,B)=1$. On the other hand, it is easy to see $\N(u) =AB$.  
Put
\[
\xi=a+bi,\qquad \zeta=d-ci.
\]
 The quaternion equations are equivalent to 
\[
\xi\zeta=u,\qquad
N(\xi)=A,\qquad N(\zeta)=B
\]
By Lemma \ref{lem:gaussian-factorization}, these have exactly four solutions. The map from $(\xi,\zeta)$ to $(a,b,c,d)$ is bijective. Hence $\# \mathcal{F}_p(X,Y,Z)=4$. 
\end{proof}

\begin{lemma}\label{lem:imprimitive-fibres}
Let $p$ be an odd prime. Suppose
\[
X^2+Y^2+Z^2=p^2
\]
and $p| \gcd(X,Y,Z)$.
Then
\[
(X,Y,Z)\in
\{(\pm p, 0,0), (0,\pm p, 0), (0,0,\pm p)\}
\]
and $\# \mathcal{F}_p(X,Y,Z)=4\left(1+\left(\frac{-1}{p}\right)\right)$ for each $(X,Y,Z)$. 
\end{lemma}

\begin{proof}
Write
\[
(X,Y,Z)=p(X_0,Y_0,Z_0).
\]
Then
\[
X_0^2+Y_0^2+Z_0^2=1,
\]
which gives the six stated possibilities. By conjugation by $h$ and $j$, it suffices to consider $(X,Y,Z)=(p,0,0)$. The equations
\[
N(\rho)=p,\qquad \mathcal T(\rho)=(p,0,0)
\]
imply
\[
a^2+b^2=p,\qquad c=d=0.
\]
Thus there are $4\left(1+\left(\frac{-1}{p}\right)\right)$ quaternions by the classical two-square formula. 
\end{proof}

\subsection{The ternary congruence conditions}\label{subsec:ternary-congruences}
Define
\begin{equation}\label{eq:Sf}
S_f(n) =\left\{ 
(X,Y,Z)\in\mathbb Z^3\colon X^2+Y^2+Z^2=n, Z\equiv0\pmod5
\right\}, 
\end{equation} 
\begin{equation}\label{eq:Sg}
S_g(n)=\left\{
(X,Y,Z)\in\mathbb Z^3\colon X^2+Y^2+Z^2=n, X+Y+2Z\equiv0\pmod5
\right\}.
\end{equation}
We set
\[
s_f(n)=\#S_f(n), \qquad s_g(n)=\#S_g(n).
\]

\begin{lemma}\label{f-bijection}
There exists a bijection 
\begin{equation}
R_f(n) \stackrel{\sim}{\longrightarrow} S_f(n)
\end{equation}
\end{lemma}

\begin{proof} 
Consider the map
\[
\Phi\colon R_f(n)\longrightarrow S_f(n),
\qquad
(x,y,z)\longmapsto(x,y,5z). 
\]
If $(x,y,z)\in  R_f(n)$, then
\[
x^2+y^2+(5z)^2=n,
\qquad 5z\equiv0\pmod5.
\]
Hence $(x,y,5z)\in S_f(n)$, so $\Phi$ is well-defined.

We then show that \(\Phi\) is a bijection.
The map is clearly injective.
To prove surjectivity, let $(X,Y,Z)\in S_f(n)$. Since $Z\equiv0\pmod5$, we can write $Z=5z$ for some $z\in\mathbb{Z}$. Setting $x=X$ and $y=Y$, we obtain
\[
f(x,y,z)
=X^2+Y^2+25\left(
\dfrac{Z}{5}
\right)^2
=X^2+Y^2+Z^2=n.
\]
Thus $(x,y,z)\in R_f(n)$ and $\Phi(x,y,z)=(X,Y,Z)$. Therefore $\Phi$ is a bijection.
\end{proof}

\begin{lemma}\label{g-bijection}
For every integer $n\geq 0$, there exists a bijection 
\begin{equation}
R_g(n) \stackrel{\sim}{\longrightarrow} S_g(n).
\end{equation}
\end{lemma}
\begin{proof}
Let $(x,y,z)\in R_g(n)$ and put
\begin{equation}
X=x+y,\qquad
Y=-x+y-z,\qquad
Z=-y-2z.
\end{equation}
By direct calculation,
\[
\begin{aligned}
X^2+Y^2+Z^2
&=(x+y)^2+(-x+y-z)^2+(-y-2z)^2\\
&=g(x,y,z)\\
&=n.
\end{aligned}
\]
Moreover,
\[
\begin{aligned}
X+Y+2Z
&=(x+y)+(-x+y-z)+2(-y-2z)\\
&=-5z.
\end{aligned}
\]
Thus
\[
X+Y+2Z\equiv0\pmod5. 
\]
Hence $(X,Y,Z)\in S_g(n)$.

Conversely, suppose that $(X,Y,Z)\in S_g(n)$. Define
\begin{equation}
\begin{aligned}
z&=-\frac{X+Y+2Z}{5},\\
y&=-Z-2z,\\
x&=X-y.
\end{aligned}
\end{equation}
Since $X+Y+2Z\equiv0\pmod5$, we have $z\in\mathbb{Z}$, and consequently $x$, $y\in\mathbb{Z}$. By construction $X=x+y$, $Y=-x+y-z$ and $Z=-y-2z$. Therefore 
\begin{equation}
g(x,y,z)=X^2+Y^2+Z^2=n,
\end{equation}
so $(x,y,z)\in R_g(n)$.
\end{proof}

\subsection{The lattices}
We define 
\begin{equation}
\Lambda_\alpha=\left\{
q_\alpha(x,y,z,w)=x+\dfrac{w}{2}+(y+\dfrac{w}{2})i+(5z+\dfrac{5w}{2})j+\dfrac{5w}{2}k \colon x,y,z,w \in \mathbb{Z}
\right\}
\end{equation}
The lattice $\Lambda_\alpha$ is contained in $\mathcal{H}$.  

Put
\[
q_0=\frac{-3+i+j+k}{2}.
\]

We define
\begin{equation}
\Lambda_\beta =\left\{ 
q_\beta(x,y,z,w)= xq_0-yq_0i+z(2j+k)+w(j-2k)\colon 
x,y,z,w\in\mathbb{Z}
\right\}.
\end{equation}
Since \(q_0\in\mathcal H\), the lattice $\Lambda_\beta$ is contained in $\mathcal{H}$.

\begin{lemma}\label{alpha-bijection}
The map
\begin{equation}
\Phi \colon \mathbb Z^4 \longrightarrow \Lambda_\alpha,  \qquad  (x,y,z,w)\longmapsto
q_\alpha(x,y,z,w)
\end{equation}
is a bijection. Moreover $\N(q_\alpha)=\alpha(x,y,z,w)$.
\end{lemma}
\begin{proof}
We can rewrite that $q_\alpha(x,y,z,w)=(x-2w)+(y-2w)i+5zj+5wh$. From the definition of $q_\alpha(x,y,z,w)$, we can find $(x,y,z,w) \in \mathbb{Z}^4$, hence $\Phi$ is surjective. 
If $\Phi(x,y,z,w)=\Phi(x^\prime,y^\prime,z^\prime,w^\prime)$ for $(x, y,z,w)\neq (x^\prime, y^\prime, z^\prime, w^\prime)$, hence by definition, we get 
\begin{equation}
(x-2w)+(y-2w)i+5zj+5wh=(x^\prime-2w^\prime)+(y^\prime-2w^\prime)i+5z^\prime j+5w^\prime h.
\end{equation}
 Since $\{ 1, i,j,h \}$ is a basis, we 
then compare the coefficients and obtain $(x,y,z,w)=(x^\prime,y^\prime,z^\prime,w^\prime)$. Hence $\Phi$ is injective. 
The direct calculation shows $\N(q_\alpha)=\alpha(x,y,z,w)$. 
\end{proof}

\begin{proposition}\label{alpha-lattice}
For $\rho=a+bi+cj+dk \in \mathcal{H}$, we have 
\begin{equation}
\rho \in \Lambda_\alpha \Longleftrightarrow c\equiv d \equiv 0\pmod 5. 
\end{equation}

\end{proposition}
\begin{proof}
If $\rho \in \Lambda_\alpha$, then we can write 
\begin{equation}
\rho=x+\dfrac{w}{2}+(y+\dfrac{w}{2})i+(5z+\dfrac{5w}{2})j+\dfrac{5w}{2}k,
\end{equation}
 hence $c\equiv d\equiv 0\pmod5$. 

Conversely, suppose these congruences hold. We set $c=\dfrac{5s}{2}$, $d=\dfrac{5t}{2}$ with $s$ and $t \in \mathbb{Z}$. 
We rewrite $\rho$ as follows 
\begin{equation}
\rho=(a-\dfrac{t}{2})+\dfrac{t}{2}+[(b-\dfrac{t}{2})+\dfrac{t}{2}]i+[5(\dfrac{s}{2}-\dfrac{t}{2})+5\dfrac{t}{2}]j+\dfrac{5t}{2}k. 
\end{equation}
Since $\rho\in \mathcal{H}$, then $2a\equiv 2b \equiv s\equiv t \pmod 2$. So $x=a-\dfrac{t}{2}$, $y=b-\dfrac{t}{2}$, $z=\dfrac{s}{2}-\dfrac{t}{2}$ and $w=t$ are all integers. Hence $\rho\in \Lambda_\alpha$. 
\end{proof}

\begin{lemma}\label{alpha-lattice-bijection}
For every integer $n\geq 0$, the map
\[
\Phi_n\colon R_\alpha(n)\longrightarrow \{
\rho \in\Lambda_\alpha\colon \N(\rho)=n
\},
\qquad
(x,y,z,w)\longmapsto q_\alpha(x,y,z,w)
\]
is a bijection.
\end{lemma}

\begin{proof}
If $(x,y,z,w)\in R_\alpha(n)$,
then
\[
\mathrm N(q_\alpha(x,y,z,w))
=\alpha(x,y,z,w)=n.
\]
Hence the map $\Phi$ in Lemma \ref{alpha-bijection} restricts to a well-defined injective map
\[
\Phi_n \colon R_\alpha(n)\longrightarrow
\{
\rho\in\Lambda_\alpha\colon \N(\rho)=n
\}.
\]

Conversely, let $\rho\in\Lambda_\alpha$ satisfy $\N (\rho)=n$. Since $\Phi$ is a bijection, there exists a unique $(x,y,z,w)\in\mathbb{Z}^4$ such that $\rho=q_\alpha(x,y,z,w)$.
Therefore
\[
n=\N (\rho)=\alpha(x,y,z,w),
\]
which shows that $(x,y,z,w)\in R_\alpha(n)$. Thus $\Phi_n$ is also surjective and hence is a bijection.
\end{proof}

\begin{lemma}\label{beta-bijection}
The map
\begin{equation}
\Phi^\prime \colon \mathbb{Z}^4 \longrightarrow \Lambda_\beta \qquad (x,y,z,w)\longmapsto
q_\beta(x,y,z,w)
\end{equation}
is bijective and $\N(q_\beta)=\beta(x,y,z,w)$.
\end{lemma}

\begin{proof}
We can rewrite $q_\beta=\dfrac{-3x+y}{2}+\dfrac{x+3y}{2}i+(\dfrac{x-y}{2}+2z+w)j+(\dfrac{x+y}{2}+z-2w)k$. 
From the definition of $q_\beta$, we can find $(x,y,z,w)\in \mathbb{Z}^4$, hence $\Phi^\prime$ is surjective. 
Now we suppose $\Phi^\prime(x, y,z,w)=\Phi^\prime(x^\prime, y^\prime, z^\prime, w^\prime)$ for $(x, y,z,w)\neq (x^\prime, y^\prime, z^\prime, w^\prime)$. 
Hence \begin{align}
\dfrac{-3x+y}{2}+\dfrac{x+3y}{2}i+(\dfrac{x-y}{2}+2z+w)j&+(\dfrac{x+y}{2}+z-2w)k= \dfrac{-3x^\prime+y^\prime}{2}+\dfrac{x^\prime+3y^\prime}{2}i\\&+
(\dfrac{x^\prime-y^\prime}{2}+2z^\prime+w^\prime)j               +(\dfrac{x^\prime+y^\prime}{2}+z^\prime-2w^\prime)k.
\end{align}
Since $\{ 1, i,j,k \}$ is a basis, we then compare the coefficients and obtain $(x,y,z,w)=(x^\prime,y^\prime,z^\prime,w^\prime)$. Hence $\Phi^\prime$ is injective. 
 We can directly get  $\N(q_\beta)=\beta(x,y,z,w)$.
\end{proof}

\begin{proposition}\label{beta-lattice}
For $\rho=a+bi+cj+dk \in \mathcal{H}$, we have 
\begin{equation}
\rho\in \Lambda_\beta \Longleftrightarrow
\begin{cases}
b-3a\equiv0\pmod5,\\
a+2c+d\equiv0\pmod5.
\end{cases}
\end{equation}
\end{proposition}
\begin{proof}
Suppose that $\rho \in\Lambda_\beta$. Then we can write 
\begin{equation}
\rho=\dfrac{-3x+y}{2}+\dfrac{x+3y}{2}i+(\dfrac{x-y}{2}+2z+w)j+(\dfrac{x+y}{2}+z-2w)k
\end{equation}
 for $(x,y,z,w)\in \mathbb{Z}^4$. So we have
\begin{align*}
b-3a&=\dfrac{x+3y}{2}-3\dfrac{-3x+y}{2}=5x \\
a+2c+d&=\dfrac{-3x+y}{2}+2(\dfrac{x-y}{2}+2z+w)+(\dfrac{x+y}{2}+z-2w)=5z. 
\end{align*}
Hence
\begin{equation}
b-3a\equiv0\pmod 5,
\qquad
a+2c+d\equiv0\pmod 5.
\end{equation}

Conversely, suppose these two congruences hold. Write $\rho=a+bi+cj+dk \in \mathcal{H}$. Then $a=\dfrac{a_0}{2}$, $b=\dfrac{b_0}{2}$, $c=\dfrac{c_0}{2}$, $d=\dfrac{d_0}{2}$ for $a_0$, $b_0$, $c_0$ and $d_0\in \mathbb{Z}$, then $a_0\equiv b_0\equiv c_0\equiv d_0 \pmod 2$ since $\rho \in \mathcal{H}$. So $b-3a$, $a+3b$, $a+2c+d$ and $b+c-2d$ are integers. 
 The first congruence gives $b\equiv3a\pmod5$. Thus
\[
a+3b\equiv a+9a=10a\equiv0\pmod5.
\]
The second congruence gives $d \equiv -a-2c\pmod 5$. Consequently, 
\begin{equation}
b+c-2d \equiv5a+5c\equiv0\pmod5.
\end{equation}
We set
 \begin{equation}
\begin{aligned}
x&=\dfrac{b-3a}{5},&
y&=\dfrac{a+3b}{5},\\
z&=\dfrac{a+2c+d}{5},&
w&=\dfrac{b+c-2d}{5}.
\end{aligned}
\end{equation}
Hence $(x,y,z,w)\in \mathbb{Z}^4$. So we can write 
\begin{equation}
\rho=a+bi+cj+dk=xq_0-yq_0i+z(2j+k)+w(j-2k).
\end{equation}
 Hence $\rho\in\Lambda_\beta$.
\end{proof}

\begin{lemma}\label{beta-lattice-bijection}
For every integer $n\geq 0$, the map
\[
\Phi_n^\prime \colon R_\beta(n)\longrightarrow
\{
\rho\in\Lambda_\beta\colon \N(\rho)=n
\},
\qquad
(x,y,z,w)\longmapsto q_\beta(x,y,z,w)
\]
is a bijection.
\end{lemma}

\begin{proof}
If $(x,y,z,w)\in R_\beta(n)$,
then
\[
\N(q_\beta(x,y,z,w)) =\beta(x,y,z,w)=n.
\]
Hence the map $\Phi^\prime$ in Lemma \ref{beta-bijection} restricts to a well-defined injective map
\[
\Phi_n^\prime \colon R_\beta(n)\longrightarrow
\{
\rho\in\Lambda_\beta\colon \N(\rho)=n
\}.
\]

Conversely, let $\rho\in\Lambda_\beta$ satisfy $\N (\rho)=n$. Since $\Phi^\prime$ is a bijection, there exists a unique $(x,y,z,w)\in\mathbb{Z}^4$ such that $\rho=q_\beta(x,y,z,w)$.
Therefore
\[
n=\N (\rho)=\beta(x,y,z,w),
\]
which shows that $(x,y,z,w)\in R_\beta(n)$.
Thus $\Phi_n^\prime$ is also surjective and hence is a bijection.
\end{proof}

\begin{lemma}\label{representation}
For every $\rho\in\mathcal{H}$, 
\begin{equation}
\rho\in\Lambda_\beta
\Longleftrightarrow
\rho\equiv q_0(s+ti)\pmod 5
\quad\text{for some } s,t\in \mathbb{F}_5.
\end{equation}
\end{lemma}

\begin{proof}
Suppose that $\rho\equiv q_0(s+ti)\pmod 5$. 
For $s,t\in\mathbb{F}_5$, we have
\[
q_0(s+ti)=\dfrac{-3s-t}{2}+\dfrac{s-3t}{2}i +\dfrac{s+t}{2}j+\dfrac{s-t}{2}k.
\]
If we write $\rho=a+bi+cj+dk$, then $b-3a\equiv 5s\equiv0\pmod5$ and $a+2c+d\equiv 0\pmod 5$. Hence, by Proposition \ref{beta-lattice}, every such $\rho$ lies in  $\Lambda_\beta$.

Conversely, suppose $\rho=a+bi+cj+dk\in\Lambda_\beta$. By Proposition \ref{beta-lattice}, we get 
\begin{equation}\label{congruence}
b\equiv3a\pmod5,
\qquad
a+2c+d\equiv0\pmod5.
\end{equation}
Set $s=c+d$ and $t=c-d$ in $\mathbb{F}_5$. Then we can write $q_0(s+ti)=-2c-d+(-c+2d)i+cj+dk$. So 
 $\rho\equiv -2c-d+(-c+2d)i+cj+dk\equiv q_0(s+ti)\pmod 5$.
\end{proof}

\subsection{Proof of the identity for $g$}\label{subsec:proof-g-beta}
For a prime $p \nmid 10$, we define
\begin{equation}\label{eq:Sg0}
S_g^{(0)}(p^2)=
\left\{
(X,Y,Z)\in S_g(p^2)\colon  (X,Y,Z)\equiv(4p,3p,4p)\pmod 5
\right\}
\end{equation} 
and $s_g^{(0)}(p^2)=\# S_g^{(0)}(p^2)$.

\begin{lemma}\label{lem:g-residue-classes}
For every prime $p\nmid10$,
\[
r_g(p^2)=4s _g^{(0)}(p^2).
\]
\end{lemma}

\begin{proof}
From Lemma \ref{g-bijection}, we have $r_g(p^2)=s_g(p^2)$. We consider $(X,Y,Z)\in S_g(p^2)$.
We have $X^2+Y^2+Z^2= p^2$ and $X+Y+2Z\equiv0\pmod5$.
Since 
\begin{equation}
Z\equiv-\dfrac{X+Y}{2}\pmod5,
\end{equation}
 we have $XY\equiv2p^2\pmod5$.
From $p\nmid 5$, we know that $X$ and $Y$ are nonzero modulo $5$.
A direct enumeration in $\mathbb F_5$ shows that
the possible residue classes are precisely
\[
\begin{array}{c|c}
\text{class}&(X,Y,Z)\pmod5\\ \hline
0&(4p,3p,4p)\\
1&(p,2p,p)\\
2&(3p,4p,4p)\\
3&(2p,p,p).
\end{array}
\]

Consider the following integral transformations
\[
\sigma\colon (X,Y,Z)\longmapsto(Y,X,Z),
\]
and
\[
\tau\colon (X,Y,Z)\longmapsto(-X,-Y,-Z).
\]
These two transformations act transitively on the four residue classes listed above. Since both transformations preserve  $X^2+Y^2+Z^2$ and $X+Y+2Z\equiv0\pmod5$, 
 they induce bijections between the corresponding subsets of $S_g(p^2)$. Consequently, all four subsets have the same cardinality, namely $s_g^{(0)}(p^2)$. This proves the result. 
\end{proof}

\begin{lemma}\label{lem:beta-inverse-integrality}
Let $p$ be a prime with $p\nmid 10$ and $(X,Y,Z)\in S_g^{(0)}(p^2)$.
Then every $\rho\in \mathcal{F}_p(X,Y,Z)$ lies in $\{\rho\in \Lambda_\beta\colon \N(\rho)=p\}$.
\end{lemma}

\begin{proof}
Since $(X,Y,Z)\equiv(4p,3p,4p)\pmod5$, we have $T(\rho) \equiv 4pi+3pj+4pk \pmod5$. Recall that
\[
\N(q_0)=3,\qquad
T(q_0)=2i-j+2k.
\]
Then it follows that
\begin{equation}\label{equality}
T(\rho) \equiv\dfrac{p}{3} T( q_0) \pmod5.
\end{equation}
Both $\rho$ and $q_0$ are invertible modulo $5$, because their norms are $p$ and $3$, respectively. We define
\[
u=q_0^{-1}\rho\pmod 5.
\]
Then $\N(u)\equiv\dfrac{p}{3}\pmod5$. In the Equation \eqref{equality}, we multiply on the left by $q_0^{-1}$ and on the right by $\bar{q}_0^{-1}$ and $u$, we obtain $\N(u)ui=\N(u)iu$. As $\N(u)\ne0$ in $\mathbb{F}_5$, we have $ui=iu$. Now we write $u=s+ti+vj+wk$ with $s$, $t$, $v$, $w\in\mathbb{F}_5$. Since $u$ commutes with $i$, we get $v=w=0$. Hence $u=s+ti$. Consequently, we get 
\[
\rho \equiv q_0(s+ti)\pmod5.
\]
By Lemma \ref{representation}, this gives $\rho\in \{\rho\in \Lambda_\beta\colon \N(\rho)=p\}$.
\end{proof}

\begin{lemma}\label{lem:beta-four-to-one}
For every prime $p\nmid10$,
\[
r_\beta(p)=4s_g^{(0)}(p^2).
\]
\end{lemma}

\begin{proof}
Suppose $(x,y,z,w) \in R_\beta(p)$. From Lemma \ref{beta-lattice-bijection}, there is a unique $\rho=q_\beta(x,y,z,w)$ with $\N(\rho)=p$.  Lemma \ref{representation} gives that $\rho\equiv q_0(s+ti) \pmod 5$ for some $s$, $t\in \mathbb{F}_5$. Taking norms, we get 
\begin{equation}
p\equiv 3(s^2+t^2) \pmod 5. 
\end{equation}
Using this, we can calculate 
\begin{equation}
T(\rho)\equiv 4pi+3pj+4pk \pmod 5. 
\end{equation}
Hence this gives a well-defined map $\mathcal{G}\colon R_\beta(p) \longrightarrow S_g^{(0)}(p^2)$. 

Let $(X,Y,Z)\in S_g^{(0)}(p^2)$. 
\begin{enumerate}
\item If $\gcd(X,Y,Z)= p$, we get $(X,Y,Z)$ must be one of the following 
\[
(\pm p,0,0),\quad(0,\pm p,0),\quad(0,0,\pm p).
\]
None of these vectors satisfies $X+Y+2Z\equiv0\pmod5$, since $p\neq5$. 
\item If $\gcd(X,Y,Z)=1$, then from Lemma \ref{lem:primitive-fibres}, we have $\# \mathcal{F}_p(X,Y,Z)=4$. By Lemma \ref{lem:beta-inverse-integrality}, all the elements in $\mathcal{F}_p(X,Y,Z)$ lie in $\{ \rho\in \Lambda_\beta, \N(\rho)=p\}$. Then $\Phi_p^{\prime -1}(\mathcal{F}_p(X,Y,Z))$ gives exactly four elements in $R_\beta(p)$ from Lemma \ref{beta-lattice-bijection}. 
\end{enumerate}
 Therefore, we have shown that the map $\mathcal{G}$ is exactly four-to-one. So we obtain the result. 
\end{proof}

\begin{lemma} \label{g-beta}
For every prime $p \nmid 10$, we have 
\[
r_g(p^2)=r_\beta(p).
\]
\end{lemma}

\begin{proof}
Combining Lemmas~\ref{lem:beta-four-to-one}
and \ref{lem:g-residue-classes}, we obtain
\[
r_g(p^2)=4s_g^{(0)}(p^2) =r_\beta(p).
\]
\end{proof}

\subsection{Proof of the identity for $f$}\label{subsec:proof-f-alpha}
Define
\begin{equation}\label{eq:Sf0}
S_f^{(0)}(p^2)= \left\{
(X,Y,Z)\in S_f(p^2)\colon  (X,Y,Z) \equiv (p,0,0) \pmod 5
\right\}
\end{equation}
and $s_f^{(0)}(p^2)=\# S_f^{(0)}(p^2)$.

\begin{lemma}\label{lem:f-four-residue-classes}
For every prime $p\nmid10$,
\[
r_f(p^2)=4 s_f^{(0)}(p^2).
\]
\end{lemma}

\begin{proof}
By Lemma \ref{f-bijection}, we have $r_f(p^2)=s_f(p^2)$. Now we consider $(X,Y,Z)\in S_f(p^2)$. Since $Z\equiv0\pmod5$,  we have 
\[
X^2+Y^2\equiv p^2\pmod5.
\]
Since $p\ne5$, we have $p^2\equiv1$ or $4\pmod5$. If neither $X$ nor $Y$ is divisible by $5$, then $X^2+Y^2\equiv 0$, $2$ or $3\pmod5$, so we get a contradiction. Hence 
\[
(X,Y, Z)\equiv
(\pm p,0, 0),\quad(0,\pm p, 0)\pmod5.
\]

The integral automorphisms
\[
(X,Y,Z)\longmapsto(-X,Y,Z),
\qquad
(X,Y,Z)\longmapsto(Y,X,Z)
\]
preserve $X^2+Y^2+Z^2$ and $Z\equiv 0 \pmod 5$.  The group generated by these two automorphisms acts transitively on the four residue classes. Hence they induce bijections between the corresponding subsets of $S_f(p^2)$. Consequently, these subsets have the same cardinality, namely $s_f^{(0)}(p^2)$. This proves the result. 
\end{proof}

\begin{lemma}\label{lem:alpha-inverse-integrality}
Let $p$ be a prime with $p\nmid 10$ and $(X,Y,Z)\in S_f^{(0)}(p^2)$. Then every $\rho \in \mathcal{F}_p(X,Y,Z)$ lies in $\{\rho \in \Lambda_\alpha, \N(\rho)=p\}$.
\end{lemma}

\begin{proof}
Let $(X,Y,Z)\in S_f^{(0)}(p^2)$. Since $(X,Y,Z)\equiv(p,0,0)\pmod5$, we have 
\begin{equation}
T(\rho) \equiv pi \pmod5. 
\end{equation}
 The quaternion $\rho$ is invertible modulo $5$, because its norm is $p$. Hence $T(\rho)\rho\equiv pi\rho \pmod 5$ and we get $\rho i\equiv i\rho\pmod 5$ after dividing $\N(\rho)$. We write $\rho=a+bi+cj+dk$. Since $\rho$ commutes with $i$ modulo $5$, we obtain $c\equiv d \equiv 0 \pmod 5$. By Proposition \ref{alpha-lattice}, this gives $\rho\in\{\rho \in \Lambda_\alpha, \N(\rho)=p\}$. 
\end{proof}

\begin{lemma}\label{lem:alpha-fibre-count}
For every prime $p\nmid10$,
\begin{equation}
r_\alpha(p)=4s_f^{(0)}(p^2) + 4\left(
\dfrac{-1}{p}
\right).
\end{equation}
\end{lemma}

\begin{proof}
Let $(x,y,z,w)\in R_\alpha(p)$. From Lemma \ref{alpha-lattice-bijection},  we have a unique $\rho\in\Lambda_\alpha$ and $\N(\rho)=p$. Proposition \ref{alpha-lattice} gives that $\rho$ commutes with $i$ modulo $5$. Then we get 
\[
T(\rho) \equiv i \rho \bar \rho=pi\pmod5.
\]
Thus we get a well-defined map $\mathcal{F}\colon R_\alpha(p)\longrightarrow S_f^{(0)}(p^2)$. 

Let $(X,Y,Z)\in S_f^{(0)}(p^2)$. We can see there is only one element $(X,Y,Z)=(p,0,0)$ with $\gcd(X,Y,Z)=p$. Every other element satisfies $\gcd(X,Y,Z)=1$. 
\begin{enumerate}
\item If $\gcd(X,Y, Z)=p$, then $\# \mathcal{F}_p(X,Y,Z)=4(1+(\dfrac{-1}{p}))$ from Lemma \ref{lem:imprimitive-fibres}. We know that Lemma \ref{lem:alpha-inverse-integrality} gives that every $\rho\in  \mathcal{F}_p(X,Y,Z)$ lies in $\{
\rho\in \Lambda_\alpha, \N(\rho)=p
\}$. Hence $\Phi_p^{-1}(\mathcal{F}_p(X,Y,Z))$ gives exactly $4(1+(\dfrac{-1}{p}))$ elements in $R_\alpha(p)$. 
\item If $\gcd(X,Y, Z)=1$, then $\# \mathcal{F}_p(X,Y,Z)=4$ from Lemma \ref{lem:primitive-fibres}. We know that Lemma \ref{lem:alpha-inverse-integrality} gives that every $\rho \in  \mathcal{F}_p(X,Y,Z)$ lies in $\{
\rho\in \Lambda_\alpha, \N(\rho)=p
\}$. Hence $\Phi_p^{-1}(\mathcal{F}_p(X,Y,Z))$ gives exactly $4$ elements in $R_\alpha(p)$ according to Lemma \ref{alpha-lattice-bijection}. 
\end{enumerate}
Consequently, we have
\begin{align*}
r_\alpha(p)
&=4\bigl(s_f^{(0)}(p^2)-1\bigr)
+4\left(1+\left(\dfrac{-1}{p}\right)\right)\\
&=
4s_f^{(0)}(p^2)+4\left(\dfrac{-1}{p}\right).
\end{align*}
\end{proof}

\begin{lemma}\label{f-alpha}
For every prime $p \nmid 10$, we have 
\[
r_f(p^2) = r_\alpha(p)-4\left(\frac{-1}{p}\right).
\]
\end{lemma}

\begin{proof}
By Lemma~\ref{lem:f-four-residue-classes} and Lemma~\ref{lem:alpha-fibre-count}, we get the result.
\end{proof}

\section{Shimura lifting}\label{sec:shimura}
Recall 
\[
\psi(z)  =\sum \limits_{n\ge1}b(n)q^n \in S_{3/2}(100,\chi_\text{triv}), 
\qquad
\phi(z)   =\sum\limits_{n\ge1}a(n)q^n \in S_2(50, \chi_\text{triv}).
\]
In particular, $b(1)=a(1)=1$.

For a positive squarefree integer $s$, we let
\begin{equation}\label{eq:chi-s}
 \xi_s(m)= \chi_\text{triv}(m) \cdot  \left(\frac{-s}{m}\right)
\end{equation}
Here we use the Kronecker-symbol convention.  

Following the classical Shimura correspondence \cite{Shimura1973}, for each positive squarefree integer $s$, we define the Shimura lift of $\psi$ by
\[
\text{Sh}_s(\psi)(z) =\sum_{n\geq 1} A_s(n)q^n, \qquad A_s(n)=\sum_{d\mid n} \xi_s(d)b\left(s\left(\frac{n}{d}\right)^2\right).
\]

\begin{theorem}\label{thm:shimura-cusp-hecke}
For every positive squarefree integer $s$, we have 
\begin{equation}
\text{Sh}_s(\psi) \in S_2(50, \chi_\text{triv}).
\end{equation}
Moreover, for every prime $p\nmid10$,
\begin{equation}
T_p\text{Sh}_s(\psi)=\lambda_p\text{Sh}_s(\psi).
\end{equation}
\end{theorem}

\begin{proof}
Since $U(100,\chi_\text{triv})=\{0\}$ from Lemma \ref{lem:no-unary}, Purkait's formulation of the Shimura
correspondence \cite[Theorem 3.1(i), (iii)]{Purkait2014} therefore gives $\text{Sh}_s(\psi) \in S_2(50, \chi_\text{triv})$.
By \cite[Theorem 1.2]{Purkait2014}, we have $\text{Sh}_s(T_{p^2}\psi)=T_p\text{Sh}_s(\psi)$.  Hence we obtain
\begin{equation}
T_p\text{Sh}_s(\psi)=\text{Sh}_s(T_{p^2}(\psi))=\text{Sh}_s(\lambda_p\psi)=\lambda_p\text{Sh}_s(\psi).
\end{equation}
\end{proof}

\begin{remark}
The theorem does not assert that $\text{Sh}_s(\psi)$ is nonzero for every squarefree $s$. Whenever the lift is nonzero, it is a simultaneous Hecke eigenform for all primes $p\nmid10$.
\end{remark}

\begin{proposition}\label{prop:hecke-eigenvalues}
For every prime $p\nmid 10$, we have $\lambda_p=a(p)$.
\end{proposition}

\begin{proof}
The coefficient of $q$ in $T_{p^2}(\psi)$ gives
\[
\lambda_p=b(p^2)+\left(\dfrac{-1}{p}\right)b(1).
\]
Since $b(1)=1$ and from the definition of $\psi$,
\begin{equation}\label{eq:lambda-representation}
\lambda_p=\dfrac{r_f(p^2)-r_g(p^2)}4+\left(\dfrac{-1}{p}\right).
\end{equation}
By Lemma \ref{f-alpha} and Lemma \ref{g-beta}, we obtain
\begin{align*}
\lambda_p &= \dfrac{r_\alpha(p)-r_\beta(p)-4(\frac{-1}{p})}{4} +\left(\dfrac{-1}{p}\right)\\
&= \dfrac{r_\alpha(p)-r_\beta(p)}4\\
&=a(p).
\end{align*}
\end{proof}

\begin{corollary}
We have
\[
\operatorname{Sh}_1(\psi)=\phi.
\]
In particular, $\phi$ is a normalized newform of weight two, level $50$, and trivial character.
\end{corollary}

\begin{proof}
Let $F=\operatorname{Sh}_1(\psi)$. Since $b(1)=1$, we have
\[
F=q+\sum_{n\geq 2}A_1(n)q^n.
\]
By Theorem \ref{thm:shimura-cusp-hecke}, we have $T_p(F)=\lambda_p F$ for $p\nmid 10$.  By Proposition \ref{prop:dim50}, the space $S_2(50,\chi_{\mathrm{triv}})$ is two-dimensional and entirely new. Let $f_1$ and $f_2$ be two distinct normalized newforms spanning this space.
Write
\[
F=c_1f_1+c_2f_2.
\]
Since $F$ is a simultaneous Hecke eigenform for all $T_p$ with $p\nmid 10$, we have
\[
c_i(a_p(f_i)-\lambda_p)=0
\]
for $i=1,2$ and every prime $p\nmid 10$. If both $c_1$ and $c_2$ were nonzero,
then $f_1$ and $f_2$ would have the same
Hecke eigenvalues at every prime $p\nmid 10$.
By strong multiplicity one, this would imply
$f_1=f_2$, a contradiction.
Thus one of $c_1,c_2$ vanishes.
Since $A_1(1)=1$, we conclude that
$F$ is a normalized newform. Now let $G$ be the other normalized newform,
so that $F$ and $G$ form a $\mathbb{C}$-basis of $S_2(50,\chi_{\mathrm{triv}})$. Since $\phi\in S_2(50,\chi_{\mathrm{triv}})$, we may write
\[
\phi=uF+vG
\]
for some $u,v\in\mathbb C$.
Comparing the Fourier coefficients of $q$ gives
$u+v=1$, and hence
\[
\phi=F+v(G-F).
\]

By Proposition \ref{prop:hecke-eigenvalues}, for every prime
$p\nmid 10$, we have $a(p)=\lambda_p=A_1(p)$. Comparing the coefficients of $q^p$ yields
\[
v\bigl(a_G(p)-a_F(p)\bigr)=0.
\]
If $v\neq 0$, then
\[
a_G(p)=a_F(p)
\qquad (p\nmid 10).
\]
By strong multiplicity one for normalized newforms, this implies $F=G$, a contradiction.
Therefore $v=0$, and hence
\[
\phi=F=\operatorname{Sh}_1(\psi).
\]
In particular, $\phi$ is a normalized newform of weight two, level $50$, and trivial character.
\end{proof}

\begin{remark}
The newform $\phi$ is associated with the elliptic curve
\[
E\colon y^2+xy+y=x^3+x^2-3x+1,
\]
which has conductor $50$ and Cremona label $50b1$. In particular, $L(\phi,s)=L(E,s)$. See also \cite[Section 8.4, Example 1]{Purkait2013}.
\end{remark}

\section{Waldspurger formula}\label{sec:waldspurger}
\subsection{The general formula}

For a nonzero integer $D$ that is a fundamental discriminant,
let $\chi_D$ denote the associated primitive quadratic Dirichlet
character and use the notation
\[
L(F,\chi_D,s)=L(F\otimes\chi_D,s).
\]

\begin{theorem}[Waldspurger]\label{thm:qin-waldspurger}
Let $\kappa\ge1$ be an integer, let $4\mid N$, and let $\chi$
be an even Dirichlet character modulo $N$. Suppose that
\[
h(z)=\sum_{n\ge1}c(n)q^n \in S_{\kappa+1/2}(N,\chi)
\]
is a simultaneous eigenform for $T_{p^2}$ for all $p\nmid N$.
Let $F\in S_{2\kappa}(N/2,\chi^2)$ be a normalized
newform with the same Hecke eigenvalues at these primes.
Write
\[
\varepsilon_\kappa=\left(\dfrac{(-1)^\kappa}{\,\cdot\,}\right),
\qquad
\eta_t=\left(\dfrac{t}{\,\cdot\,}\right).
\]
If $t$ and $t_0$ are positive squarefree integers satisfying
\[
t/t_0\in(\mathbb{Q}_\ell^\times)^2
\qquad\text{for every prime }\ell\mid N,
\]
then
\begin{equation}
c(t_0)^2\,t^{\kappa-1/2}
 L(F,\varepsilon_\kappa^{-1}\chi^{-1}\eta_t,\kappa)
 =  c(t)^2\,t_0^{\kappa-1/2}  \chi(t_0/t)  L(F,\varepsilon_\kappa^{-1}\chi^{-1}\eta_{t_0},\kappa).
\end{equation}
\end{theorem}

\begin{proof}
This is the formula in \cite[Theorem~3.19]{Qin2022}; see also \cite[Theorem~2.9]{LuQin2025}.
\end{proof}

\begin{corollary}\label{cor:waldspurger-ours}
If $t,t_0$ are positive squarefree integers with
$t/t_0\in(\mathbb Q_\ell^\times)^2$ for $\ell=2,5$, then
\begin{equation}\label{eq:waldspurger-ours}
b(t_0)^2\sqrt t\,
 L(\phi, \chi_{D_t},1)
=
b(t)^2\sqrt{t_0}\,
 L(\phi, \chi_{D_{t_0}},1),
\end{equation}
where $D_u$ denotes the fundamental discriminant of
$\mathbb Q(\sqrt{-u})$. In particular, for even squarefree $u$,
we have $D_u=-4u$.
\end{corollary}

\begin{proof}
Apply Theorem~\ref{thm:qin-waldspurger} with $\kappa=1$, $N=100$, $\chi=\chi_\text{triv}$, $h=\psi$, $F=\phi$.
Since $\varepsilon_1=\bigl(\frac{-1}{\cdot}\bigr)$ and the associated primitive quadratic character is $\chi_{D_t}$, the required equality follows.
\end{proof}

\subsection{Local square classes at $2$ and $5$}
We consider $t=20n+2$ and we assume that $t$ is squarefree. It follows that $v_2(t)=1$, $5\nmid t$, and $D_t=-4t$ is a
negative fundamental discriminant.
Obviously, there is a unique $t_0\in\{2,22,42,62\}$ such that $t\equiv t_0\pmod{80}$.
\begin{lemma}\label{lem:local-square-classes}
We have 
\[
t/t_0\in(\mathbb Q_2^\times)^2
\quad\text{and}\quad
t/t_0\in(\mathbb Q_5^\times)^2.
\]
\end{lemma}

\begin{proof}
Write $t=2u$ and $t_0=2u_0$, with $u,u_0$ odd.
Since $t\equiv t_0\pmod{80}$,
we have $u\equiv u_0\pmod8$; therefore
\[
\frac{t}{t_0}=\frac{u}{u_0}\in1+8\mathbb Z_2.
\]
Every element of $1+8\mathbb Z_2$ is a square in
$\mathbb Q_2^\times$. Also $t\equiv t_0\equiv2\pmod5$,
so $t/t_0\in1+5\mathbb Z_5$, which consists of squares by
Hensel's lemma. 
\end{proof}

\begin{lemma}\label{lem:four-reference-coefficients}
For $t_0\in\{2,22,42,62\}$, the representation numbers are
\[
\begin{array}{c|rrrr}
t_0&2&22&42&62\\ \hline
r_f(t_0)&4&0&16&16\\
r_g(t_0)&2&6&8&20\\
b(t_0)&\frac12&-\frac32&2&-1
\end{array}
\]
\end{lemma}

\begin{proof}
The result follows from direct calculation. 
\end{proof}

\begin{theorem}\label{thm:explicit-waldspurger}
Let $t$ be a positive squarefree integer satisfying $t\equiv2\pmod{20}$.
Then
\begin{equation}\label{eq:explicit-L-value}
L(\phi, \chi_{-4t},1)
=\dfrac{b(t)^2}{\sqrt t}
\begin{cases}
\displaystyle 4\sqrt2
L(\phi, \chi_{-8},1),
 &t\equiv2\pmod{80},\\[7pt]
\displaystyle \frac{4\sqrt{22}}{9}
L(\phi,\chi_{-88},1),
 &t\equiv22\pmod{80},\\[7pt]
\displaystyle \frac{\sqrt{42}}{4}
L(\phi,\chi_{-168},1),
 &t\equiv42\pmod{80},\\[7pt]
\displaystyle \sqrt{62}
L(\phi,\chi_{-248},1),
 &t\equiv62\pmod{80}.
\end{cases}
\end{equation}
\end{theorem}

\begin{proof}
Choose $t_0$ as in Lemma~\ref{lem:local-square-classes}.
By Corollary \ref{cor:waldspurger-ours},
\[
b(t_0)^2\sqrt t\,L(\phi,\chi_{-4t},1)
=b(t)^2\sqrt{t_0}\,L(\phi,\chi_{-4t_0},1).
\]
Since $b(t_0)\ne0$, 
this gives
\[
L(\phi,\chi_{-4t},1) =\dfrac{b(t)^2}{\sqrt t} \dfrac{\sqrt{t_0}}{b(t_0)^2} L(\phi,\chi_{-4t_0},1).
\]
Substitute the four values of $b(t_0)$ from Lemma \ref{lem:four-reference-coefficients}.
\end{proof}

\begin{lemma}\label{waldspurgerbound}
For $u\in\{2,22,42,62\}$, one has
\[
L(\phi\otimes\chi_{-4u},1)\neq 0.
\]
Moreover, there is an effectively computable
constant $C_W>0$ such that
\[
|b(t)|^2
\leq C_W\sqrt{t}\,
\bigl|L(\phi\otimes\chi_{-4t},1)\bigr|
\]
for every positive squarefree $t\equiv 2\pmod{20}$.
\end{lemma}

\begin{proof}
Since $\phi$ is a normalized newform of weight $2$, level $50$, and trivial character, hypothesis (H1) follows from \cite[Theorem~6]{Purkait2013}. Moreover, since $v_2(50)=1$, \cite[Corollary 4.1]{Purkait2013} shows that
the local representation $\rho_2$ is not supercuspidal. Hence (H2) also holds. By Waldspurger's theorem \cite[Theorem~7]{Purkait2013}, the form $\psi$ is a finite linear combination of forms $f(c,A_\phi)$, where $A_\phi$ is a $\mathbb{C}$-valued function on $\{n\in N, n \text{ is squarefree} \}$. Consequently, for every positive squarefree integer $u$, its Fourier coefficient $b(u)$ is of the form
\[
b(u)=A_\phi(u)u^{1/4}
\sum_{j=1}^r\lambda_j
\prod_p c_p^{(j)}(u).
\]
See more details in \cite[Section 5]{Purkait2013}. For even squarefree $u$, the character $\chi_\text{triv}^{-1}\varepsilon_1^{-1}\chi_u$ is precisely $\chi_{-4u}$. Since the epsilon factor of a quadratic character equals $1$, we have
\[
A_\phi(u)^2=L(\phi,\chi_{-4u},1).
\]
It follows that 
\[
L(\phi,\chi_{-4u},1)=0
\quad\Longrightarrow\quad b(u)=0.
\]
On the other hand, from Lemma \ref{lem:four-reference-coefficients}, we find that $b(t_0)\neq 0$. Therefore $L(\phi,\chi_{-4u},1)\ne0$, for $u\in\{2,22,42,62\}$. By the nonvanishing established above, choosing the unique $u\in\{2,22,42,62\}$ such that $t\equiv u\pmod{80}$, we have 
\[
|b(t)|^2 =\dfrac{|b(u)|^2} {\sqrt u\,|L(\phi\otimes\chi_{-4u},1)|} \sqrt t\,|L(\phi, \chi_{-4t},1)|.
\]
Define
\[
C_W=\max_{u\in\{2,22,42,62\}}\dfrac{|b(u)|^2} {\sqrt u\,|L(\phi,\chi_{-4u},1)|}.
\]
Then $C_W>0$ and $|b(t)|^2\le C_W\sqrt t\, |L(\phi,  \chi_{-4t},1)|$.

Finally, the four fixed central values defining $C_W$ can be evaluated with certified error bounds, for instance by modular symbols or Mellin integrals. Their proven nonvanishing ensures that the required divisions and approximations are effective. Hence $C_W$ is effectively computable.
\end{proof}

\section{Application to Sun's conjecture}\label{sec:application}

In this section, we apply the results established in the preceding sections to Sun's conjecture. We begin by reducing the original representation problem to that of a positive definite ternary quadratic form. This reduction allows us to use the arithmetic properties of ternary quadratic forms developed above.
We have the following elementary equivalence. 

\begin{proposition}\label{prop:Sun-equivalence}
For every integer $n\ge0$,  we have
\begin{equation}\label{eq:Sun-equivalence}
 n=P(x,y,z)\text{ for some }x,y,z\in\mathbb{Z}  \quad\Longleftrightarrow\quad
 20n+2=u^2+v^2+25w^2\text{ for some }u,v,w\in\mathbb{Z}.
\end{equation}
\end{proposition}

\begin{proof}
Completing squares gives
\[
 40n+4=2(10x+1)^2+(10y+1)^2+(10z+1)^2.
\]
Put $u=10x+1$, $v=5(y+z)+1$ and $w=y-z$. Then
$10y+1=v+5w$, $10z+1=v-5w$, and therefore
\[
 40n+4=2u^2+(v+5w)^2+(v-5w)^2
        =2(u^2+v^2+25w^2).
\]
Conversely, suppose $u^2+v^2+25w^2=20n+2$.
Modulo $5$ one has $u^2+v^2\equiv2$, which forces
$u^2\equiv v^2\equiv1\pmod5$. Changing signs of $u$ and $v$
if necessary, arrange $u\equiv v\equiv1\pmod5$.
Modulo $4$, precisely two of $u,v,w$ are odd. Since $u$ and $v$
occur symmetrically, exchange them if necessary to make $u$ odd.
Then $u\equiv1\pmod{10}$; and since $v$ and $w$ have opposite
parities, both $v\pm5w$ are odd and congruent to $1$ modulo $5$,
thus congruent to $1$ modulo $10$.
Consequently
\[
 x=\dfrac{u-1}{10},\qquad
 y=\dfrac{v+5w-1}{10},\qquad
 z=\dfrac{v-5w-1}{10}
\]
are integers. Reversing the square completion proves $n=P(x,y,z)$.
\end{proof}

Hence the representation problem  reduces to
\[
 r_f(20n+2)>0,\qquad
 f(x,y,z)=x^2+y^2+25z^2.
\]

\subsection{The genus average and three squares}

\begin{proposition}\label{prop:genus-three-square}
For every positive integer $n\equiv2, 3\pmod5$,
\begin{equation}\label{eq:genus-three-square}
r_f(n)+4r_g(n)=r_{Q_0}(n), 
\end{equation}
where $Q_0=x^2+y^2+z^2$.
\end{proposition}

\begin{proof}
We know that $5\Theta=\theta_f+4\theta_g$ from the definition. Hence $5a_n(\Theta)=r_f(n)+4r_g(n)$, where $\Theta=\sum \limits_{n\geq 0} a_n(\Theta) q^n$. 
By the Minkowski--Siegel formula in \cite[Theorem 2.1]{BerkovichPatane2015}, we have 
\[
a_n(\Theta)= \dfrac{2\pi\sqrt{n}}{5} \prod_{\ell}\alpha_\ell(n,f),
\]
where $\alpha_\ell(n,f)= \lim_{k\to\infty}\ell^{-2k} \#\left\{
\bm{x} \in(\mathbb{Z}/\ell^k\mathbb{Z})^3\colon f(\bm{x})\equiv n\pmod{\ell^k}
\right\}$.

Since $\text{gen}(Q_0)=\{Q_0\}$, from Minkowski--Siegel formula again, we have 
\[
r_{Q_0}(n) = 2\pi\sqrt{n} \prod_{\ell}\alpha_\ell(n,Q_0),
\]
where $\alpha_\ell(n,Q_0)= \lim_{k\to\infty}\ell^{-2k} \#\left\{
\bm{x} \in(\mathbb{Z}/\ell^k\mathbb{Z})^3\colon Q_0(\bm{x})\equiv n \pmod{\ell^k}
\right\}$.

For every prime $\ell\ne5$, the change of variables $(x,y,z)\mapsto(x,y,5z)$ is invertible over $\mathbb Z_\ell$. Hence $f$ and $Q_0$ are equivalent over $\mathbb{Z}_\ell$. Then we have 
\[
\alpha_\ell(n,f)=\alpha_\ell(n,Q_0).
\]
At $\ell=5$, the assumption $n\equiv2\pmod5$ implies that
\[
x^2+y^2\equiv2\pmod5. 
\]
This has precisely four solutions for $(x,y)$, namely $(\pm1,\pm1)$. Since $z$ is arbitrary, there are $4\cdot5=20$ solutions in $\mathbb{F}_5^3$.
Similarly, the congruence $x^2+y^2+z^2\equiv2\pmod5$ has exactly $20$ solutions in $\mathbb{F}_5^3$.
Every solution in either case is nonsingular modulo $5$. Hence each solution modulo $5$ lifts to exactly $5^{2(k-1)}$ solutions modulo $5^k$. Therefore,
\[
\alpha_5(n,f)=\alpha_5(n,Q_0)=\dfrac{20}{25}=\dfrac{4}{5}.
\]
Hence we have $5a_n(\Theta)=r_{Q_0}(n)$. So we get $r_f(n)+4r_g(n)=r_{Q_0}(n)$. 

The same argument applies when $n\equiv3\pmod5$. Indeed, both $f$ and $Q_0$ have exactly $20$ nonsingular solutions modulo $5$, giving $\alpha_5(n,f)=\alpha_5(n,Q_0)=\dfrac{4}{5}$.
Thus the identity holds.

\end{proof}

Let $H(M)$ denote the Hurwitz class number of discriminant $-M$, normalized so that $H(|D|)=h(D)$ for every negative fundamental discriminant $D<-4$. See more details in \cite[Section 3]{Li2023}.

\begin{corollary}\label{cor:representation-Hurwitz}
For every positive integer $n\equiv2\pmod{20}$,
\begin{equation}\label{eq:representation-decomposition}
r_f(n)=\dfrac{12}{5}H(4n)+\dfrac{16}{5}b(n).
\end{equation}
If $n$ is squarefree, then $-4n$ is fundamental and
\begin{equation}\label{eq:representation-squarefree}
r_f(n)=\dfrac{12}{5}h(-4n)+\dfrac{16}{5}b(n).
\end{equation}
\end{corollary}

\begin{proof}
By the definitions of $\Theta$ and $\psi$, we have $\theta_f=\Theta+\frac{16}{5}\psi$. Taking the coefficient of $q^n$, we obtain
\[
r_f(n)=a_n(\Theta)+\dfrac{16}{5}b(n).
\]
Since $n\equiv2\pmod{20}$, Proposition 7.2 and Gauss's three-square formula \cite[\text{Theorem 1.4}]{Mortenson2017} give
\[
a_n(\Theta) =\dfrac{1}{5}r_{Q_0}(n)=\dfrac{12}{5}H(4n),
\]
which proves the first identity.

If $n$ is squarefree, then $n\equiv2\pmod4$ implies that $-4n$ is a fundamental discriminant. Since $-4n<-4$, we have $H(4n)=h(-4n)$, giving the second identity. 
\end{proof}

\subsection{When $t$ is not squarefree}

In this section, we consider positive integers $t\equiv 2\pmod{20}$ that are not squarefree. 

\begin{theorem}\label{thm:nonsquarefree}
Let $t\equiv2\pmod{20}$ be a positive integer.
If $t$ is not squarefree, then
\[
r_f(t)>0.
\]
\end{theorem}

\begin{proof}
Write $t=sk^2$, where $s$ is squarefree and $k>1$.
Since $t\equiv2\pmod{20}$, we have
\[
s\equiv2\pmod4,\qquad (k,10)=1,
\]
and $s\equiv2$ or $3\pmod5$. 
Since $s\equiv2\pmod4$, Legendre's three-square theorem gives $r_{Q_0}(s)>0$. By Proposition \ref{prop:genus-three-square} and $s\equiv2,3\pmod5$, we have $r_f(s)+4r_g(s)=r_{Q_0}(s)>0$. If $r_f(s)>0$, we can choose $(x_1, y_1, z_1)$ satisfying $x_1^2+y_1^2+25z_1^2=s$. Hence we have $x_1^2k^2+y_1^2k^2+25z_1^2k^2=sk^2$, which shows that $(x_1k, y_1k, z_1k)$ belongs to $R_f(t)$. So $r_f(t)>0$. 

Suppose now that $r_f(s)=0$.
Recall that
\[
\theta_f=\Theta+\frac{16}{5}\psi.
\]
Comparing the Fourier coefficient of $q^s$, we obtain $r_f(s)=a_s(\Theta)+\dfrac{16}{5}b(s)$. Hence $b(s)=-\dfrac{5}{16}a_s(\Theta)$. 
Moreover, the proof of Proposition \ref{prop:genus-three-square} gives $a_s(\Theta)=\dfrac{1}{5} r_{Q_0}(s)>0$.
Choose a prime $p| k$.
Since $(k,10)=1$, we have $p\nmid10$. The Hecke relations $T_{p^2}$ for $\Theta$ and $\psi$ give
\[
a_{sp^2}(\Theta)
=(p+1-\xi_s(p))a_s(\Theta)
\]
and
\[
b(sp^2)=(a(p)-\xi_s(p))b(s).
\]
Consequently
\begin{align*}
r_f(sp^2)
&=a_{sp^2}(\Theta)+\dfrac{16}{5}b(sp^2)\\
&=\bigl(p+1-\xi_s(p)\bigr)a_s(\Theta)\\
&\quad-\bigl(a(p)-\xi_s(p)\bigr)a_s(\Theta)\\
&=\dfrac{1}{5}\bigl(p+1-a(p)\bigr) r_{Q_0}(s).
\end{align*}
By Deligne's bound, $|a(p)|\le2\sqrt p$.
Therefore,
\[
p+1-a(p)\ge p+1-2\sqrt p
=(\sqrt p-1)^2>0.
\]
We conclude that $r_f(sp^2)>0$. Then we can choose $(x_1^\prime, y_1^\prime, z_1^\prime)$ in $R_f(sp^2)$. Hence the vector $(x_1^\prime k/p, y_1^\prime k/p, z_1^\prime k/p)$ belongs to $R_f(t)$. So $r_f(t)>0$. 
\end{proof}

\begin{remark}
The proof follows the Hecke-theoretic method of Lu and Qin \cite{LuQin2017}. Their stated theorem assumes coprimality with the level, whereas the integers considered here satisfy \((t,100)=2\); the preceding argument treats this case directly.
\end{remark}

\subsection{The squarefree case: an effective bound}
We now consider positive squarefree integers $t\equiv2\pmod{20}$.

\begin{theorem}\label{thm:effective-squarefree}
There exists an effectively computable positive integer $M_0$ such that
\[
  r_f(t)>0
\]
for every squarefree integer $t\geq M_0$ satisfying
$t\equiv2\pmod{20}$.
\end{theorem}
\begin{proof}
Let $t$ be squarefree and satisfy $t\equiv2\pmod{20}$.
By Corollary~7.3,
\begin{equation}\label{eq:effective-representation-identity}
  r_f(t)
  =
  \frac{12}{5}h(-4t)+\frac{16}{5}b(t).
\end{equation}
We estimate the two terms on the right with effectively computable
constants.

For the required class-number estimate, we use \cite[Theorem~1.3(ii)]{BanksBS} directly.
Let $c_{\mathrm B}>0$ be the effectively computable constant in that theorem for the base field $\mathbb Q$.
The imaginary quadratic field $K=\mathbb Q(\sqrt{-t})$
has
\[
  d_K=4t,\qquad   [K:\mathbb Q]=2,\qquad
  R_K=1,\qquad
  h_K=h(-4t).
\]
Consequently,
\[
  h(-4t)
  \geq
  c_{\mathrm B}
  \dfrac{\sqrt{4t}}{(2\pi)^2\log(4t)}.
\]
Putting $c_h=\dfrac{c_{\mathrm B}}{2\pi^2}$, we obtain the effective lower bound
\begin{equation}\label{eq:effective-class-number}
  h(-4t)
  \geq
  c_h\frac{\sqrt t}{\log(4t)}.
\end{equation}

We next estimate $b(t)$.
By Lemma~6.6, there is an effectively computable constant $C_W>0$
such that
\begin{equation}\label{eq:effective-waldspurger-bound}
  |b(t)|^2
  \leq
  C_W\sqrt t\,
  \bigl|L(\phi,\chi_{-4t},1)\bigr|.
\end{equation}
Since $\phi$ is a fixed normalized newform of weight $2$, level $50$, and trivial nebentypus, the subconvexity estimate of Blomer and Harcos \cite[Theorem 2]{BlomerHarcos}, together with its addendum \cite{BlomerHarcosAddendum}, gives, for every $\varepsilon>0$,
\begin{equation}\label{eq:effective-subconvexity}
  \bigl|L(\phi,\chi_{-4t},1)\bigr|
  \leq
  C_{\phi,\varepsilon}\,t^{3/8+\varepsilon},
\end{equation}
where $C_{\phi,\varepsilon}>0$ is effectively computable.

Taking $\varepsilon=1/16$ in
\eqref{eq:effective-subconvexity}, we obtain
\[
  \bigl|L(\phi,\chi_{-4t},1)\bigr|
  \leq
  C_L\,t^{7/16},
  \qquad
  C_L=C_{\phi,1/16}.
\]
Combining this with Equation \eqref{eq:effective-waldspurger-bound} gives
\[
  |b(t)|^2\leq C_WC_L\,t^{15/16}.
\]
Thus, with
\[
  C_b=(C_WC_L)^{1/2},
\]
we have
\begin{equation}\label{eq:effective-fourier-bound}
  |b(t)|\leq C_b\,t^{15/32}.
\end{equation}
Both $C_b$ and $c_h$ are effectively computable.
It follows from Equation \eqref{eq:effective-class-number} and Equation \eqref{eq:effective-fourier-bound} that
\[
  \frac{4|b(t)|}{3h(-4t)}
  \leq
  \frac{4C_b}{3c_h}\,
  t^{-1/32}\log(4t).
\]
Choose a positive integer $A$ effectively such that
\[
  A\geq
  \max\left\{1,\frac{4C_b}{3c_h}\right\}
\]
and put $M_0=4(128A)^{64}$. For $x\geq1$, we have the elementary inequality
\[
  \log x\leq64x^{1/64}. 
\]
 So whenever $t\geq M_0$,
\[
  A\,t^{-1/32}\log(4t)
  \leq
  64A\,4^{1/64}t^{-1/64}
  \leq\frac12.
\]
Therefore
\[
  \dfrac{4|b(t)|}{3h(-4t)}\leq\dfrac{1}{2}
  \qquad (t\geq M_0).
\]
Finally, Equation \eqref{eq:effective-representation-identity} gives
\begin{align*}
  r_f(t)   &\geq
  \dfrac{12}{5}h(-4t)-\dfrac{16}{5}|b(t)|\\
  &=
  \dfrac{12}{5}h(-4t)
  \left(
    1-\dfrac{4|b(t)|}{3h(-4t)}
  \right)\\
  &\geq
  \dfrac65h(-4t)\\
  &\geq
  \dfrac{6c_h}{5}\dfrac{\sqrt t}{\log(4t)}
  >0.
\end{align*}

\end{proof}

\begin{corollary}\label{cor:Sun}
There exists an effectively computable positive integer $N_0>0$ such that every integer $n\ge N_0$ can be represented in the form
\[
n=P(x,y,z)
=x(5x+1)+\frac{y(5y+1)}{2}
+\frac{z(5z+1)}{2}
\]
for some $x,y,z\in\mathbb{Z}$.
\end{corollary}

\begin{proof}
Let $M_0$ be the effectively computable positive integer in Theorem \ref{thm:effective-squarefree}. We put
\[
N_0=\max\left\{ 1,\left\lceil\dfrac{M_0-2}{20}\right\rceil
\right\}.
\]
Let $n\geq N_0$ and set $t=20n+2$. Then $t\geq M_0$ and $t\equiv2\pmod{20}$.
If $t$ is squarefree, Theorem \ref{thm:effective-squarefree} gives $r_f(t)>0$.
If $t$ is not squarefree, the same conclusion follows from Theorem \ref{thm:nonsquarefree}. Proposition \ref{prop:Sun-equivalence} therefore gives an integral representation of $n$ by $P$. Since $M_0$ is effectively computable, so is $N_0$.
\end{proof}

\newpage

\bibliographystyle{amsplain}
\bibliography{sun_conjecture_1}

\end{document}